\documentclass[12pt, reqno]{amsart}
\usepackage{amsmath, amsthm, amscd, amsfonts, amssymb, graphicx, xcolor}
\usepackage[bookmarksnumbered, colorlinks, plainpages]{hyperref}
\usepackage{amsmath}
\usepackage{cases}
\usepackage{amssymb}
\usepackage{mathrsfs}
\usepackage{euscript}
\usepackage{dsfont}
\usepackage{graphicx}
\usepackage{float}
\usepackage{graphicx}
\usepackage{amssymb}
\newcommand{\pp}{\prec\mkern-5mu\prec}

\allowdisplaybreaks[4]
\newtheorem{theorem}{Theorem}[section]
\newtheorem{thmletter}{Theorem}

\newtheorem{lemma}{Lemma}[section]
\newtheorem{proposition}[theorem]{Proposition}

\theoremstyle{definition}

\theoremstyle{remark}
\newtheorem{remark}[theorem]{Remark}
\numberwithin{equation}{section}

\begin{document}
	\setcounter{page}{1}
	

	
	\title[]{Smooth numbers in simultaneous arithmetic progressions}

	\author[]{Zhiyuan Yang}
	
	\address{School of Mathematics,  Shandong University, Jinan 250100, Shandong, China}
	\email{\tt zhiyuan.yang@mail.sdu.edu.cn}

	
	

	
	
	\begin{abstract}
		We prove a mean value theorem on the distribution of smooth numbers in two simultaneous
		arithmetic progressions, with generalizing a result of Drappeau. Let $\varepsilon>0$ and $(\log x)^C<y<x^{1/C}$ for a sufficiently large $C$. Our result implies that $y$-smooth numbers are well-distributed in the residue classes of a fixed integer $a$ on average over moduli $\leq x^{5/9-\varepsilon}$ and $\leq x^{153/256-\varepsilon}$ for each $|a|\leq x$, which improves Drappeau's range $\leq x^{6/11-\varepsilon}$. As an application, we show that a positive proportion of $y$-smooth numbers $n\leq x$ satisfy $\min(P^+(n-1), P^+(n+1))>x^{407/720}$.
		\newline
		\noindent 
	\end{abstract} \maketitle
	
	
	\section{Introduction}
	\subsection{Smooth numbers in arithmetic progressions}
	Define $\pi(x) := \#\{p \leq x\}$. For $(a,q)=1$, let $\pi(x;q,a) := \#\{p \leq x : p \equiv a\mkern-7mu\pmod{q}\}$. The
	celebrated Bombieri-Vinogradov theorem [3, 20] asserts that for $B = B(A)$ large enough
	in terms of $A$, one has
	\begin{align*}
		\sum_{\substack{q\leq x^{1/2}/(\log x)^B}}\max_{(a,q)=1}\Bigg|\pi(x;q,a)-\frac{\pi(x)}{\varphi(q)}\Bigg|\ll_A\frac{x}{(\log x)^A}.
	\end{align*}
	For an integer $n > 1$, let $P^+(n)$ be the largest prime factor of $n$ with the convention that $P^+(1)=1$.
	Replacing primes by other numbers, Wolke [21] proved a theorem of Bombieri-Vinogradov type for sifted numbers, and also announced that
	\begin{align*}
		\sum_{\substack{q\leq x^{1/2}/(\log x)^B}}\max_{(a,q)=1}\Bigg|\Psi(x,y;q,a)-\frac{\Psi_q(x,y)}{\varphi(q)}\Bigg|\ll_A\frac{x}{(\log x)^A},
	\end{align*}
	uniformly for $x\geq y\geq 2$, where
	\begin{align*}
		\Psi(x,y;q,a):=\sum_{\substack{n\leq x,P^+(n)\leq y\\n\equiv a\mkern-15mu\pmod{q}}}1,\quad \Psi_q(x,y):=\sum_{\substack{n\leq x,P^+(n)\leq y\\(n,q)=1}}1.
	\end{align*}
	We also define $\Psi(x,y):=\Psi_1(x,y)$.
	In 1991, Fouvry and Tenenbaum [7] proved the above result with a different method. Later, Harper gave a better bound in [10]. In order to obtain a level of the distribution of $>1/2$ for smooth numbers for some small $y$, Fouvry and Tenenbaum [8] in 1996 also proved that
	\begin{thmletter}[Fouvry-Tenenbaum]
		For any $\varepsilon>0$, there exists $\delta>0$ such that the following holds. For any $A>0$ and any $1\leq y\leq x^\delta$, for $a\in\mathbb{Z}\backslash\{0\}$, one has
		\begin{align*}
			&\sum_{\substack{q\leq x^{3/5-\varepsilon}\\(q,a)=1}}\Bigg|\Psi(x,y;q,a)-\frac{\Psi_q(x,y)}{\varphi(q)}\Bigg|\ll\frac{x}{(\log x)^A},\quad \text{for}\ |a|\leq x^\delta,\tag{1.1}\\&\sum_{\substack{q\leq x^{6/11-\varepsilon}\\(q,a)=1}}\Bigg|\Psi(x,y;q,a)-\frac{\Psi_q(x,y)}{\varphi(q)}\Bigg|\ll\frac{x}{(\log x)^A},\quad \text{for}\ |a|\leq x,\tag{1.2}
		\end{align*}
		and
		for any well-factorable function $\lambda(q)$ of level $x^{5/9-\varepsilon}$, \begin{align*}
			\sum_{(q,a)=1}\lambda(q)\Bigg(\Psi(x,y;q,a)-\frac{\Psi_q(x,y)}{\varphi(q)}\Bigg)\ll\frac{x}{(\log x)^A},\quad \text{for}\ |a|\leq x.\tag{1.3}
		\end{align*}
	\end{thmletter}
	To understand (1.3), it is necessary to introduce what it means to be well-factorable. A complex sequence $(\lambda(q))_{q\leq Q}$ is said to be
	well-factorable of level $Q$ iff for any $Q_1$, $Q_2\geq 1$ with $Q_1Q_2=Q$, there exist 1-bounded
	complex sequences $(\alpha_{q_1})$, $(\beta_{q_2})$ supported on $q_i \leq Q_i$, such that for all $q$,
	\begin{align*}
		\lambda(q)=\sum_{q_1q_2=q}\alpha_{q_1}\beta_{q_2}.
	\end{align*}
	This weight arises in the Rosser-Iwaniec sieve [12] and has found numerous applications.
	In 2015, Drappeau [4] proved the above results (1.1) and (1.2) with a different method and also obtained a better bound, which allows for the larger range $(\log x)^C< y<x^{1/C}$. Very recently, for (1.1), Pascadi improved the level to $x^{66/107-\varepsilon}$ and $x^{5/8-\varepsilon}$ respectively in [18] and [19]. In this paper, we improve the level in (1.2) to $x^{5/9-\varepsilon}$.
	
	Furthermore, we wish to extend the above result to cover more values of $a$. A typical example is under two congruence classes. 
	In 2021, Assing, Blomer and Li [1, Theorem 2.1] showed that
	\begin{align*}
		\sum_{\substack{q\leq Q\\(q,a_1)=1}}\Bigg(\sum_{\substack{p\leq x\\p\equiv a_1\mkern-15mu\pmod{q}\\p\equiv d_0\mkern-15mu\pmod{d}}}1-\frac{\pi(x)}{\varphi(dq)}\Bigg)\ll\frac{x}{(\log x)^A},
	\end{align*}
	where $Q\leq x^{1/2+\delta}$, $|a_1|\leq x^{1+\delta}$, $d\leq (\log x)^C$ for some suitable $\delta$ and $C$.
	Later, in 2023, Blomer, Grimmelt, Li and Rydin Myerson [2, Propositions 6.1 and 7.1] showed the following results, for which we give some simpler versions.
	\begin{thmletter}[Blomer-Grimmelt-Li-Rydin Myerson]
		Suppose $|\lambda_d|\ll\tau(d)^C$ for some ﬁxed $C > 0$. Then there exist some constants $\varpi,D > 0$ such that the following holds. Let $x \geq 2$, $a_2\in\mathbb{Z}\backslash\{0\}$ such that for
		\begin{align*}
			Q\leq x^{1/2+\varpi},\quad |a_2|\leq x^{1+\varpi},\quad (\log x)^D\leq y\leq x,
		\end{align*}
		we have
		\begin{align*}
			&\sum_{d\leq x^{\varpi}}\lambda_d\max_{(a_1,d)=1}\Bigg|\sum_{\substack{q\leq Q\\(q,a_2d)=1}}\Bigg(\sum_{\substack{p\leq x\\p\equiv a_1\mkern-15mu\pmod{d}\\p\equiv a_2\mkern-15mu\pmod{q}}}1-\frac{\pi(x)}{\varphi(dq)}\Bigg)\Bigg|\ll\frac{x}{(\log x)^A},\\
			&\sum_{d\leq x^{\varpi}}\lambda_d\max_{(a_1,d)=1}\Bigg|\sum_{\substack{q\leq Q\\(q,a_2d)=1}}\Bigg(\sum_{\substack{n\leq x,P^+(n)\leq y\\n\equiv a_1\mkern-15mu\pmod{d}\\n\equiv a_2\mkern-15mu\pmod{q}}}1-\frac{1}{\varphi(dq)}\sum_{\substack{n\leq x,P^+(n)\leq y\\(n,qd)=1}}1\Bigg)\Bigg|\ll\frac{\Psi(x,y)}{(\log x)^A}.
		\end{align*}
	\end{thmletter}
	As an application, they studied
	\begin{align*}
		p=x_1^2+x_2^2+1\ \text{and}\ n=x_1^2+x_2^2+m,
	\end{align*}
	where $x_1$, $x_2$ are almost primes and $m$ is $y$-smooth. 
	Unfortunately, the above result has certain limited scope. For example, the result in Theorem B is extremely inconvenient for sieve problems. Very recently, Zheng [23] studied \begin{align*}
		\sum_{d}\sum_{\substack{q\\(q,d)=1}}\gamma_d\lambda_q\Bigg(\sum_{\substack{p\leq x\\p\equiv a_1\mkern-15mu\pmod{q}\\p\equiv a_2\mkern-15mu\pmod{d}}}1-\frac{\pi(x)}{\varphi(dq)}\Bigg)\ll\frac{x}{(\log x)^A},
	\end{align*}
	where $|a_1|\leq (\log x)^B$, $|a_2|\leq x$ and $a_1\neq a_2$, $\boldsymbol{\gamma}=(\gamma_d)$ and $\boldsymbol{\lambda}=(\lambda_q)$ are divisor-bounded sequences, $\boldsymbol{\lambda}$ is well-factorable (arising in the Rosser-Iwaniec sieve).
	
	Motivated by the above results, we present our following result.
	\begin{theorem}
		For any $\varepsilon>0$, there exist $\delta$, $C>0$ such that the following holds. Let $x\geq 2$ and 
		$|a_2|\leq x^{1+\delta}$. Then for any $y\in[(\log x)^C, x^{1/C}]$ and $A>0$, one has
		\begin{align*}
			\sum_{\substack{q\leq Q}}\max_{(a_1,q)=1}\sum_{\substack{r\leq R\\(r,a_2q)=1}}\Bigg|\sum_{\substack{n\leq x,P^+(n)\leq y\\n\equiv a_1\mkern-15mu\pmod{q}\\n\equiv a_2\mkern-15mu\pmod{r}}}1-\frac{1}{\varphi(qr)}\sum_{\substack{n\leq x,P^+(n)\leq y\\(n,qr)=1}}1\Bigg|\ll\frac{\Psi(x,y)}{(\log x)^A},
		\end{align*}
		where $Q\leq x^{1/10}$ and $R=x^{5/9-\varepsilon}Q^{-14/9}$.
	\end{theorem}
	
	We are also able to replace the above upper bound by $\Psi(x,y)(H(u)^{-\delta}(\log x)^{-A}+y^{-\delta})$, where $u=\log x/\log y$, $H(u)=\exp(u/(\log(u+1))^2)$.
	
	\begin{remark}
		If we relax the condition $|a_2|\leq x$ to $|a_2|\leq x^{o(1)}$, using the methods of Drappeau and Pascadi in [4, 19] together with generalizations of several tools, the author can also prove the following result.
		For any $\varepsilon>0$, there exists $\delta$, $C>0$ such that the following holds. Let $x\geq 2$ and 
		$|a_2|\leq x^{\delta}$. Then for any $y\in[(\log x)^C, x^{1/C}]$ and $A\geq 0$, one has
		\begin{align*}
			\sum_{\substack{q\leq Q}}\max_{(a_1,q)=1}\sum_{\substack{r\leq R\\(r,a_2q)=1}}\Bigg|\sum_{\substack{n\leq x,P^+(n)\leq y\\n\equiv a_1\mkern-15mu\pmod{q}\\n\equiv a_2\mkern-15mu\pmod{r}}}1-\frac{1}{\varphi(qr)}\sum_{\substack{n\leq x,P^+(n)\leq y\\(n,qr)=1}}1\Bigg|\ll\frac{\Psi(x,y)}{(\log x)^A},
		\end{align*}
		where $Q\leq x^{25/449}=x^{1/18+1/8082}$ and $R=x^{5/8-\varepsilon}Q^{-649/200}$. However, Theorem 1.1 suffices for most problems, if one wants to avoid extremely complicated calculations.
	\end{remark}
	
	Furthermore, setting $q=1$, we also prove the level can be as large as $x^{(5-\theta)/8-o(1)}$ by using Lichtman’s generalized result for Kloosterman sums.
	\begin{theorem}
		For any $\varepsilon>0$, there exist $\delta$, $C>0$ such that the following holds. Let $x\geq 2$ and 
		$|a|\leq x^{1+\delta}$. Then for any $y\in[(\log x)^C, x^{1/C}]$ and $A>0$, one has
		\begin{align*}
			\sum_{\substack{r\leq x^{153/256-\varepsilon}\\(r,a)=1}}\Bigg|\sum_{\substack{n\leq x,P^+(n)\leq y\\n\equiv a\mkern-15mu\pmod{r}}}1-\frac{1}{\varphi(r)}\sum_{\substack{n\leq x,P^+(n)\leq y\\(n,r)=1}}1\Bigg|\ll\frac{\Psi(x,y)}{(\log x)^A}.
		\end{align*}
	\end{theorem}
	Note that $153/256\approx 0.598>5/9\approx 0.556$. Pascadi’s refined tools [19, Propositions 3.4 and 3.7] yield no benefit here.
	
	\subsection{A direct application of Theorem 1.1}
	We give the following definition. An integer $n$ is said to be $strong$ if $n-1$ has a large prime factor $q$, $n+1$ has a large prime factor $q'$. More precisely, we call an integer $n$\ $(\theta_1,\theta_2)$-$strong$ if
	\begin{align*}
		P^+(n-1)>x^{\theta_1}, \quad P^+(n+1)>x^{\theta_2}.
	\end{align*}
	This was first studied by Gordon [9]. 
	In 2015, Meng [17] showed that
	the proportion of $(1/2,1/2)$-$strong$ primes is at least $1/8$.
	By appealing to the Brun-Titchmarsh theorems on average, Xi [22] showed that  a positive proportion of primes are
	$(1/2+1/36,1/2+1/36)$-$strong$. Motivated by their works, we ask how $strong$ smooth numbers can be. 
	
	Our result is as follows.
	\begin{theorem}
		There exists $C>0$ such that for any $y\in[(\log x)^C, x^{1/C}]$, we have
		\begin{align*}
			\#\{2\leq n\leq x:P^+(n)\leq y,P^+(n-1)>x^{t}, P^+(n+1)>x^{t}\}\gg \Psi(x,y),
		\end{align*}
		for any $t<407/720=1/2+1/16+1/360$.
	\end{theorem}
	
	Combining Remark 1.2, we can improve the constant $407/720\approx 0.565$ to $0.566$.

	\section{Lemmas}
	In this section, we prepare the necessary lemmas.
	\begin{lemma}
		For any pair $a$, $b$ of non-zero coprime integers, one has
		\begin{align*}
			\frac{\overline{a}}{b}+\frac{\overline{b}}{a}\equiv\frac{1}{ab}\mkern-15mu\pmod{1}.
		\end{align*}
	\end{lemma}
	\begin{proof}
		We can see it in [6, Lemma 6]. 
	\end{proof}
	\begin{lemma}
		For all $r$, $q$, $h\geq 1$ and $\chi\mkern-7mu\pmod{r}$, we have
		\begin{align*}
			\Bigg|\sum_{\substack{0<b<q\\(b,q)=1}}\chi(b)e\left(-\frac{bh}{q}\right)\Bigg|\ll\mathrm{cond}(\chi)^{1/2}\sum_{d|(h,q)}d.
		\end{align*}
	\end{lemma}
	\begin{proof}
		It follows easily from Gauss sums and Ramanujan sums.
	\end{proof}
	
	\begin{lemma}
		For $k$, $a$, $b$, $c$, $l\in \mathbb{N}$ and $D\geq 1$, we have
		\begin{align*}
			&\sum_{\substack{d\leq D\\(d,c)=1}}e\left(b\frac{\overline{d}}{c}\right)\ll\log(D+1)(b,c)^{1/2}\tau(c)c^{1/2}+\frac{(b,c)}{c}D,
		\end{align*}
		\begin{align*}
			\sum_{\substack{d\leq D\\(d,ck)=1\\d\equiv0\mkern-15mu\pmod{l}}}\frac{d}{\varphi(d)}e\left(b\frac{\overline{d}}{c}\right)\ll\left(\log(D+1)^2(b,c)^{1/2}\tau(c)c^{1/2}+\frac{(b,c)\log(l+1)}{lc}D\right)2^{\omega(k)}.
		\end{align*}
		For $0<A_2-A_1\leq c$, for all integers $\lambda$ and $\Lambda\geq 1$, we also have
		\begin{align*}
			\sum_{\substack{A_1<d\leq A_2\\(d,ac)=1\\d\equiv\lambda\mkern-15mu\pmod{\Lambda}}}e\left(b\frac{\overline{d}}{c}\right)\ll(b,c)^{1/2}c^{1/2}\tau(ac)\log(2c).
		\end{align*}
	\end{lemma}
	\begin{proof}
		We can see it in [5, Lemma 4] and [4, Lemma 3].
	\end{proof}

	\begin{lemma}[Large sieve inequality]
		For all integers $Q$, $M$, $N\geq 1$ and any sequence of complex numbers $(a_n)_{M<n\leq M+N}$, we have
		\begin{align*}
			\sum_{q\leq Q}\frac{q}{\varphi(q)}\sum_{\substack{\chi\mkern-15mu\pmod{q}\\ \chi\ \text{primitive}}}\Bigg|\sum_{M<n\leq M+N}a_n\chi(n)\Bigg|^2\leq (N+Q^2-1)\sum_{M<n\leq M+N}|a_n|^2.
		\end{align*}
	\end{lemma}
	\begin{proof}
		We can see it in [13, Theorem 7.13].
	\end{proof}
	
	\begin{lemma}
		There exist constants $\eta$, $c>0$ such that for $(\log x)^c\leq y\leq x$ and $Q\leq x$, for any $A>0$, we have\begin{align*}
			\sum_{q\leq Q}\frac{\tau(q)}{\varphi(q)}\sum_{\substack{\chi\mkern-15mu\pmod{q}\\1<\mathrm{cond}(\chi)\leq x^{\eta}}}\Bigg|\sum_{n\leq x,P^+(n)\leq y}\chi(n)\Bigg|\ll\frac{\Psi(x,y)}{(\log x)^A}.
		\end{align*}
	\end{lemma}
	\begin{proof}
		For the estimation, we proceed identically as in [10,
		after (2.3)] with $q_0=1$, $y_2=x$.
	\end{proof}
	
	\begin{lemma}[Truncated Poisson summation]
		Let $x\gg 1$ and $1\ll N,Q\ll x^{O(1)}$, $a\in \mathbb{Z}$, $q\in \mathbb{Z}_+$ with $q\asymp Q$, and $\Phi:(0,\infty)\rightarrow\mathbb{C}$ be a smooth function, $\Phi(t)$ supported in $t\asymp 1$, with $\Phi^{(k)}\ll_k1$ for $k\geq 0$. Then for any $A$, $\delta>0$ and $H:=x^{\delta}N^{-1}Q$, one has\begin{align*}
			\sum_{n\equiv a\mkern-15mu\pmod{q}}\Phi\left(\frac{n}{N}\right)=\frac{N}{q}\sum_{|h|\leq H}\hat{\Phi}\left(h\frac{N}{q}\right)e\left(-\frac{ah}{q}\right)+O_{A,\delta}(x^{-A}),
		\end{align*}
		where\begin{align*}
			\hat{\Phi}(s):=\int\Phi(t)e(st)\mathrm{d}t.
		\end{align*}
	\end{lemma}
	\begin{proof}
		This is easy from the Poisson summation formula.
	\end{proof}

	\begin{lemma}
		Given $y\geq 2$, $N_1, N_2\geq 1$, for every integer $n>yN_1N_2$ satisfying $P^+(n)\leq y$, there exists a unique factorization $n=n_0n_1n_2$ such that
		\begin{align*}
			N_2<n_2\leq N_2P^-(n_2), \ P^+(n_2)\leq y,
		\end{align*}
		\begin{align*}
			N_1<n_1\leq N_1P^-(n_1), \ P^+(n_1)\leq P^-(n_2),
		\end{align*}
		\begin{align*}
			P^+(n_0)\leq P^-(n_1).
		\end{align*}
	\end{lemma}
	\begin{proof}
		We can see it in [8, Lemma 3.2].
	\end{proof}
	
	\subsection*{Notation.}
	We use $\varepsilon$ to denote a sufficiently small positive number, and the value of $\varepsilon$ may change from statement to statement. Let $P^+(n)$ and $P^-(n)$ denote the largest and smallest prime factors of $n$ respectively with the convention $P^+(1)=1$ and $P^-(1)=\infty$. We use the standard asymptotic notation $f = o(g)$, $f = O(g)$, $f \ll g$, $f \asymp g$, indicating dependencies of implicit constants through subscripts. Define
	\begin{align*}
		A\pp B
	\end{align*}
	to mean $A\ll_{\varepsilon} x^{O(\varepsilon)}B$. Statements like  $A\ll B^{O(\varepsilon)}C$, $A\ll B^{-O(\varepsilon)}C$ should respectively be read as $A\ll B^{c_1\varepsilon} C$, $A\ll B^{-c_2\varepsilon}C$  for some sufficiently large absolute constants $c_1$, $c_2>0$. Statements like  $AB^{-O(\varepsilon)}\ll C$, $AB^{O(\varepsilon)}\ll C$ should respectively be read as $A\ll B^{O(\varepsilon)}  C$, $A\ll B^{-O(\varepsilon)}C$. We use $n \sim N$ for the range of integers $n \in (N, 2N]$. Define $e(z):=\exp(2\pi iz)$. We denote Euler's function by $\varphi(n)$. Let $\chi$ denote a Dirichlet character.
	By $(m_1,m_2)$ we denote the greatest common divisor of $m_1$ and $m_2$.
	
	\section{Proof of Theorem 1.1}
	
	In this section, we employ the proof framework of Drappeau [4]. We mainly use Linnik's dispersion method.
	
	For all $r\in\mathbb{N}$, $k\mkern-7mu\pmod{r}$, we define
	\begin{align*}
		\omega_\varepsilon(k;r):=\sum_{\substack{\chi\ \text{primitive}\\ \text{cond}(\chi)\leq x^{\varepsilon}\\\text{cond}(\chi)|r}}\chi(k).
	\end{align*}
	In fact, we have
	\begin{align*}
		\#\{\chi\ \text{primitive}: \text{cond}(\chi)\leq x^{\varepsilon}\}\leq x^{2\varepsilon}.
	\end{align*}
	For any $a_1$ and $a_2$, define
	\begin{align*}
		&E(q,r;a_1,a_2):=\sum_{\substack{n\leq x,P^+(n)\leq y\\n\equiv a_1\mkern-15mu\pmod{q}\\n\equiv a_2\mkern-15mu\pmod{r}}}1-\frac{1}{\varphi(qr)}\sum_{\substack{n\leq x,P^+(n)\leq y\\(n,qr)=1}}1,
		\\&E_\varepsilon (q,r;a_1,a_2):=\sum_{\substack{n\leq x,P^+(n)\leq y\\n\equiv a_1\mkern-15mu\pmod{q}\\n\equiv a_2\mkern-15mu\pmod{r}}}1-\frac{1}{\varphi(qr)}\sum_{\substack{n\leq x,P^+(n)\leq y\\(n,qr)=1}}\omega_\varepsilon(n\overline{a_1};q)\omega_\varepsilon(n\overline{a_2};r).
	\end{align*}
	Then we have
	\begin{align*}
		\sum_{\substack{q\leq Q}}&\max_{(a_1,q)=1}\sum_{\substack{r\leq R\\(r,a_2q)=1}}\Big|E(q,r;a_1,a_2)\Big|\leq \sum_{\substack{q\leq Q}}\max_{(q,a_1)=1}\sum_{\substack{r\leq R\\(r,a_2q)=1}}\Big|E_\varepsilon(q,r;a_1,a_2)\Big|\\&+\sum_{\substack{q\leq Q\\(q,a_1)=1}}\max_{(a_1,q)=1}\sum_{\substack{r\leq R\\(r,a_2q)=1}}\frac{1}{\varphi(qr)}\Bigg|\sum_{\substack{n\leq x,P^+(n)\leq y\\(n,qr)=1}}\mathop{\sum_{\substack{\chi_1\ \text{primitive}\\ \text{cond}(\chi_1)\leq x^{\varepsilon}\\\text{cond}(\chi_1)|q}}\sum_{\substack{\chi_2\ \text{primitive}\\ \text{cond}(\chi_2)\leq x^{\varepsilon}\\\text{cond}(\chi_2)|r}}}_{\substack{\text{cond}(\chi_1)\text{cond}(\chi_2)>1}}\chi_1\chi_2(n)\overline{\chi_1(a_1)\chi_2(a_2)} \Bigg|
	\end{align*}
	and by Lemma 2.5,
	\begin{align*}
		&\mathop{\sum_{\substack{q\leq Q}}\max_{(a_1,q)=1}\sum_{\substack{r\leq R\\(r,a_2q)=1}}}_{}\frac{1}{\varphi(qr)}\Bigg|\sum_{\substack{n\leq x,P^+(n)\leq y\\(n,qr)=1}}\mathop{\sum_{\substack{\chi_1\ \text{primitive}\\ \text{cond}(\chi_1)\leq x^{\varepsilon}\\\text{cond}(\chi_1)|q}}\sum_{\substack{\chi_2\ \text{primitive}\\ \text{cond}(\chi_2)\leq x^{\varepsilon}\\\text{cond}(\chi_2)|r}}}_{\substack{\text{cond}(\chi_1)\text{cond}(\chi_2)>1}}\chi_1\chi_2(n)\overline{\chi_1(a_1)\chi_2(a_2)} \Bigg|\\&\leq\sum_{\substack{s\leq QR}}\sum_{\substack{s=qr\\(r,a_2q)=1}}\frac{1}{\varphi(qr)}\mathop{\sum_{\substack{\chi_1\ \text{primitive}\\ \text{cond}(\chi_1)\leq x^{\varepsilon}\\\text{cond}(\chi_1)|q}
			}\sum_{\substack{\chi_2\ \text{primitive}\\ \text{cond}(\chi_2)\leq x^{\varepsilon}\\\text{cond}(\chi_2)|r}}}_{\text{cond}(\chi_1)\text{cond}(\chi_2)>1}\Bigg|\sum_{\substack{n\leq x,P^+(n)\leq y\\(n,qr)=1}}\chi_1\chi_2(n) \Bigg|\\&\leq\sum_{s\leq QR}\frac{\tau(s)}{\varphi(s)}\sum_{\substack{\chi\mkern-15mu\pmod{s}\\ 1<\text{cond}(\chi)\leq x^{2\varepsilon}}}\Bigg|\sum_{\substack{n\leq x,P^+(n)\leq y}}\chi(n) \Bigg|\\&\ll\frac{\Psi(x,y)}{(\log x)^{A}}.
	\end{align*} 
	Since $\Psi(x,(\log x)^C)=x^{1-1/C+o(1)}$,
	we only need to prove that for some $\delta>0$, one has
	\begin{align*}
		\sum_{\substack{q\leq Q}}\max_{(a_1,q)=1}\sum_{\substack{r\leq R\\(r,a_2q)=1}}\Big|E_\varepsilon(q,r;a_1,a_2)\Big|\ll x^{1-\delta}.
	\end{align*}
	For each $q$, we suppose the maximum is attained at $a_1=a_{1q}$.
	Following the proof of [4, Proposition 2], it suffices to prove that for some $\delta>0$, one has
	\begin{align*}
		&\Delta=\Delta(Q,R,M,N,L)=\\&\sum_{\substack{q\sim Q}}\sum_{\substack{r\sim R\\(r,a_2q)=1}}\Bigg|\mathop{\sum_m\sum_n\sum_l}_{\substack{mnl\equiv a_{1q}\mkern-15mu\pmod{q}\\mnl\equiv a_2\mkern-15mu\pmod{r}}}\alpha_m\beta_n\gamma_l-\frac{1}{\varphi(qr)}\mathop{\sum_m\sum_n\sum_l}_{(mnl,qr)=1}\alpha_m\beta_n\gamma_l\omega_\varepsilon(mnl\overline{a_{1q}},q)\omega_\varepsilon(mnl\overline{a_2},r)\Bigg|\\&\ll x^{1-\delta},\tag{3.1}
	\end{align*}
	where $m\sim M$, $n\sim N$, $l\sim L$ and $MNL=x$, $|\alpha_m|\leq 1$, $|\beta_n|\leq 1$ and $|\gamma_l|\leq 1$.
	We use Linnik’s dispersion method (see [16]) to deal with $\Delta(Q,R,M,N,L)$. By Cauchy-Schwarz inequality, we have
	\begin{align*}
		&\Delta^2\leq MQR\\& \times\sum_{\substack{q\sim Q}}\sum_{\substack{r\sim R\\(r,a_2q)=1}}\sum_{\substack{m\sim M\\(m,qr)=1}}\Bigg|\sum_{\substack{k\equiv a_{1q}\overline{m}\mkern-15mu\pmod{q}\\k\equiv a_2\overline{m}\mkern-15mu\pmod{r}}}u_k-\frac{1}{\varphi(qr)}\sum_{\substack{(k,qr)=1}}u_k\omega_\varepsilon(mk\overline{a_{1q}},q)\omega_\varepsilon(mk\overline{a_2},r)\Bigg|^2 ,
	\end{align*}
	where
	\begin{align*}
		u_k=\sum_{\substack{nl=k\\n\sim N, l\sim L}}\beta_n\gamma_l.
	\end{align*}
	There exists a smooth function $f$ supported on $[M/2, 3M]$, satisfying\begin{align*}
		&\left\{ \begin{aligned}
			&f(m)=1, \quad&& M\leq m\leq 2M;\\&f(m)\geq0,\quad&& m\in \mathbb{R};\\&f^{(j)}(m)\ll_j M^{-j},\quad&& j\geq 0.
		\end{aligned}\right.
	\end{align*}
	Then
	\begin{align*}
		\Delta^2\leq MQR(\mathscr{S}_1-2Re\mathscr{S}_2+\mathscr{S}_3),\tag{3.2}
	\end{align*}
	where
	\begin{align*}
		&\mathscr{S}_1=\sum_{\substack{q\sim Q}}\sum_{\substack{r\sim R\\(r,a_2q)=1}}\sum_{(m,qr)=1}f(m)\Bigg|\sum_{\substack{k\equiv a_{1q}\overline{m}\mkern-15mu\pmod{q}\\k\equiv a_2\overline{m}\mkern-15mu\pmod{r}}}u_k\Bigg|^2,\\&\mathscr{S}_2=\sum_{\substack{q\sim Q}}\sum_{\substack{r\sim R\\(r,a_2q)=1}}\frac{1}{\varphi(qr)}\sum_{(m,qr)=1}f(m)\\&\qquad\times\sum_{\substack{(k_1,qr)=1}}\sum_{\substack{k_2\equiv a_{1q}\overline{m}\mkern-15mu\pmod{q}\\k_2\equiv a_2\overline{m}\mkern-15mu\pmod{r}}}u_{k_1}\overline{u_{k_2}}\omega_\varepsilon(mk_1\overline{a_{1q}},q)\omega_\varepsilon(mk_1\overline{a_2},r),\\&\mathscr{S}_3=\sum_{\substack{q\sim Q}}\sum_{\substack{r\sim R\\(r,a_2q)=1}}\frac{1}{(\varphi(qr))^2}\sum_{(m,qr)=1}f(m)\Bigg|\sum_{\substack{(k,qr)=1}}u_k\omega_\varepsilon(mk\overline{a_{1q}},q)\omega_\varepsilon(mk\overline{a_2},r)\Bigg|^2.
	\end{align*}
	
	In what follows, the estimates for $\mathscr{S}_2$ and $\mathscr{S}_3$ are relatively simple; we shall focus on estimating $\mathscr{S}_1$.

	\subsection{Estimate of $\mathscr{S}_3$}
	Changing the order of summation, we get
	\begin{align*}
		\mathscr{S}_3=&\sum_{\substack{q\sim Q}}\sum_{\substack{r\sim R\\(r,a_2q)=1}}\frac{1}{(\varphi(qr))^2}\sum_{\substack{\chi_1\ \text{primitive}\\\text{cond}(\chi_1)<x^\varepsilon\\\text{cond}(\chi_1)|q}}\sum_{\substack{\chi_2\ \text{primitive}\\\text{cond}(\chi_2)<x^\varepsilon\\\text{cond}(\chi_2)|q}}\overline{\chi_1}\chi_2(a_{1q})\sum_{\substack{\chi_3\ \text{primitive}\\\text{cond}(\chi_3)<x^\varepsilon\\\text{cond}(\chi_3)|r}}\sum_{\substack{\chi_4\ \text{primitive}\\\text{cond}(\chi_4)<x^\varepsilon\\\text{cond}(\chi_4)|r }}\overline{\chi_3}\chi_4(a_2)\\&\times\mathop{\sum_{k_1}\sum_{k_2}}_{(k_1k_2,qr)=1}\chi_1\chi_3(k_1)u_{k_1}\overline{\chi_2\chi_4(k_2)u_{k_2}}\sum_{(m,qr)=1}f(m)\chi_1\chi_3(m)\overline{\chi_2\chi_4(m)}.
	\end{align*}
	For $0<b<qr$, $(b,qr)=1$ and $m\equiv b\mkern-7mu\pmod{qr}$, we have
	\begin{align*}
		&\sum_{m\equiv b\mkern-15mu\pmod{qr}}f(m)\chi_1\chi_3(m)\overline{\chi_2\chi_4(m)}=\chi_1\chi_3(b)\overline{\chi_2\chi_4(b)}\sum_{m\equiv b\mkern-15mu\pmod{qr}}f(m)\\&=\frac{\chi_1\chi_3(b)\overline{\chi_2\chi_4(b)}}{qr}\sum_{|h|\leq H}\hat{f}\left(\frac{h}{qr}\right)e\left(-\frac{bh}{qr}\right)+O\left(\frac{1}{qr}\right),
	\end{align*}
	where $H=x^\varepsilon QRM^{-1}$. Then
	\begin{align*}
		\sum_{(m,qr)=1}f(m)\chi_1\chi_3\overline{\chi_2\chi_4(m)}=\sum_{\substack{0<b<qr\\(b,qr)=1}}\frac{\chi_1\chi_3\overline{\chi_2\chi_4(b)}}{qr}\sum_{|h|\leq H}\hat{f}\left(\frac{h}{qr}\right)e\left(-\frac{bh}{qr}\right)+O\left(1\right).
	\end{align*}
	By Lemma 2.2, we have \begin{align*}
		\frac{1}{qr}\sum_{0<|h|\leq H}\hat{f}\left(\frac{h}{qr}\right)\sum_{\substack{0<b<qr\\(b,qr)=1}}\chi_1\chi_3\overline{\chi_2\chi_4}(b)e\left(-\frac{bh}{qr}\right)\ll\frac{x^{2\varepsilon}}{qr}\sum_{d|qr}d\sum_{0<|h|\leq H/d}\Big|\hat{f}\left(\frac{dh}{qr}\right)\Big|\ll x^{3.1\varepsilon},
	\end{align*}
	since $|\hat{f}(dh/(qr))|\ll qr/(dh)$.
	For $h=0$, we have
	\begin{align*}
		\overline{\chi_1}\chi_2(a_{1q})\overline{\chi_3}\chi_4(a_2)\sum_{\substack{0<b<qr\\(b,qr)=1}}\chi_1\chi_3\overline{\chi_2\chi_4}(b)=\sum_{\substack{0<b<qr\\(b,qr)=1}}\chi_1\chi_3\overline{\chi_2\chi_4}(\xi b)=\sum_{\substack{0<b<qr\\(b,qr)=1}}\chi_1\chi_3\overline{\chi_2\chi_4}(b),
	\end{align*}
	where $\xi\equiv r\overline{a_{1q}r}^{(q)}+q\overline{a_2q}^{(r)}\mkern-7mu\pmod{qr}$ ($\overline{a}^{(q)}$ denotes an inverse of $a$ modulo $q$). Then we get
	\begin{align*}
		\mathscr{S}_3=\hat{f}(0)X_3+R_3,\tag{3.3}
	\end{align*}
	where
	\begin{align*}
		&X_3=\sum_{\substack{q\sim Q}}\sum_{\substack{r\sim R\\(r,a_2q)=1}}\frac{1}{qr}\sum_{\substack{0<b<qr\\(b,qr)=1}}\Bigg|\frac{1}{\varphi(qr)}\sum_{(k,qr)=1}u_k\omega_\varepsilon(kb,q)\omega_\varepsilon(kb,r)\Bigg|^2,\\&R_3\ll x^{11.1\varepsilon}\sum_{q\sim Q}\sum_{r\sim R}\frac{1}{(\varphi(qr))^2}\Big(\sum_{k}|u_k|\Big)^2\ll x^{12\varepsilon}K^2Q^{-1}R^{-1}.
	\end{align*}

	\subsection{Estimate of $\mathscr{S}_2$}Changing the order of summation, we get
	\begin{align*}
		\mathscr{S}_2=&\sum_{\substack{\chi_1\ \text{primitive}\\\text{cond}(\chi_1)<x^\varepsilon}}\sum_{\substack{\chi_2\ \text{primitive}\\\text{cond}(\chi_2)<x^\varepsilon}}\sum_{\substack{q\sim Q\\\text{cond}(\chi_1)|q}}\sum_{\substack{r\sim R\\(r,a_2q)=1\\\text{cond}(\chi_2)|r}}\frac{1}{\varphi(qr)}\sum_{(m,qr)=1}f(m)\\&\quad\times\sum_{\substack{(k_1,qr)=1}}\sum_{\substack{k_2\equiv a_{1q}\overline{m}\mkern-15mu\pmod{q}\\k_2\equiv a_2\overline{m}\mkern-15mu\pmod{r}}}u_{k_1}\overline{u_{k_2}}\chi_1(mk_1\overline{a_{1q}})\chi_2(mk_1\overline{a_2})\\=&\sum_{\substack{\chi_1\ \text{primitive}\\\text{cond}(\chi_1)<x^\varepsilon}}\sum_{\substack{\chi_2\ \text{primitive}\\\text{cond}(\chi_2)<x^\varepsilon}}\sum_{\substack{q\sim Q\\(q,a_{1q})=1\\\text{cond}(\chi_1)|q}}\sum_{\substack{r\sim R\\(r,a_2q)=1\\\text{cond}(\chi_2)|r}}\frac{1}{\varphi(qr)}\sum_{(m,qr)=1}f(m)\\&\quad\times\sum_{\substack{(k_1,qr)=1}}\sum_{\substack{k_2\equiv a_{1q}\overline{m}\mkern-15mu\pmod{q}\\k_2\equiv a_2\overline{m}\mkern-15mu\pmod{r}}}u_{k_1}\overline{u_{k_2}}\chi_1(k_1\overline{k_2})\chi_2(k_1\overline{k_2})\\=&\sum_{\substack{\chi_1\ \text{primitive}\\\text{cond}(\chi_1)<x^\varepsilon}}\sum_{\substack{\chi_2\ \text{primitive}\\\text{cond}(\chi_2)<x^\varepsilon}}\sum_{k_1}\chi_1\chi_2(k_1)u_{k_1}\sum_{k_2}\overline{\chi_1\chi_2(k_2)u_{k_2}}\\&\quad\times\sum_{\substack{q\sim Q\\(q,k_1k_2)=1\\\text{cond}(\chi_1)|q}}\sum_{\substack{r\sim R\\(r,a_2qk_1k_2)=1\\\text{cond}(\chi_2)|r}}\frac{1}{\varphi(qr)}\sum_{\substack{m\equiv a_{1q}\overline{k_1}\mkern-15mu\pmod{q}\\m\equiv a_2\overline{k_1}\mkern-15mu\pmod{r}}}f(m).
	\end{align*}
	Let $\varrho\mkern-7mu\pmod{qr}$ be a common solution of
	\begin{align*}
		&\left\{ \begin{aligned}
			&m\equiv a_{1q}\overline{k_1}\mkern-15mu\pmod{q},\\&m\equiv a_2\overline{k_1}\mkern-15mu\pmod{r}.
		\end{aligned}\right.
	\end{align*}
	We can write $\varrho\equiv a_{1q}r\overline{k_1r}^{(q)}+a_2q\overline{k_1q}^{(r)}\mkern-7mu\pmod{qr}$. Then we have
	\begin{align*}
		\sum_{\substack{m\equiv \varrho\mkern-15mu\pmod{qr}}}f(m)&=\frac{1}{qr}\sum_{|h|\leq H}\hat{f}\left(\frac{h}{qr}\right)e\left(-\frac{\varrho h}{qr}\right)+O\left(x^{-100}\right)\\&=\frac{1}{qr}\sum_{|h|\leq H}\hat{f}\left(\frac{h}{qr}\right)e\left(-\frac{a_{1q}h\overline{k_1r}}{q}-\frac{a_2h\overline{k_1q}}{r}\right)+O\left(x^{-100}\right).
	\end{align*}
	The contribution of $O(x^{-100})$ is negligible. Since \begin{align*}
		\frac{\overline{k_1q}}{r}\equiv-\frac{\overline{r}}{k_1q}+\frac{1}{k_1qr}\mkern-15mu\pmod{1},
	\end{align*}
	we get
	\begin{align*}
		e\left(-\frac{a_{1q}h\overline{k_1r}}{q}-\frac{a_2h\overline{k_1q}}{r}\right)&=e\left(-a_{1q}h\frac{\overline{k_1r}}{q}+a_2h\frac{\overline{r}}{k_1q}-\frac{a_2h}{k_1qr}\right)\\&=e\left(\varrho_1\frac{\overline{r}}{k_1q}-\frac{a_2h}{k_1qr}\right),
	\end{align*}
	where $\varrho_1\equiv a_2h-a_{1q}hk_1\overline{k_1}^{(q)}\mkern-7mu\pmod{k_1q}$. 
	Then we have
	\begin{align*}
		\mathscr{S}_2=\hat{f}(0)X_2+R_2+O(1),\tag{3.4}
	\end{align*}
	where
	\begin{align*}
		&X_2=\sum_{\substack{q\sim Q}}\sum_{\substack{r\sim R\\(r,a_2q)=1}}\frac{1}{qr\varphi(qr)}\mathop{\sum_{k_1}\sum_{k_2}}_{(k_1k_2,qr)=1}u_{k_1}\overline{u_{k_2}}\omega_\varepsilon(k_1\overline{k_2},q)\omega_\varepsilon(k_1\overline{k_2},r),
		\\&R_2=\sum_{\substack{\chi_1\ \text{primitive}\\\text{cond}(\chi_1)<x^\varepsilon}}\sum_{\substack{\chi_2\ \text{primitive}\\\text{cond}(\chi_2)<x^\varepsilon}}\sum_{k_1}\chi_1\chi_2(k_1)u_{k_1}\sum_{k_2}\overline{\chi_1\chi_2(k_2)u_{k_2}}\\&\qquad\times\sum_{\substack{q\sim Q\\(q,k_1k_2)=1\\\text{cond}(\chi_1)|q}}\sum_{\substack{r\sim R\\(r,a_2qk_1k_2)=1\\\text{cond}(\chi_2)|r}}\frac{1}{qr\varphi(qr)}\sum_{0<|h|\leq H}\hat{f}\left(\frac{h}{qr}\right)e\left(\varrho_1\frac{\overline{r}}{k_1q}-\frac{a_2h}{k_1qr}\right).
	\end{align*}
	Now we begin to estimate $R_2$. We have
	\begin{align*}
		R_2\ll&\sum_{\substack{\chi_1\ \text{primitive}\\\text{cond}(\chi_1)<x^\varepsilon}}\sum_{\substack{\chi_2\ \text{primitive}\\\text{cond}(\chi_2)<x^\varepsilon}}\sum_{k_1}\sum_{k_2}|u_{k_1}u_{k_2}|\sum_{\substack{q\sim Q\\(q,k_1k_2)=1\\\text{cond}(\chi_1)|q}}\frac{1}{q\varphi(q)}\sum_{0<|h|\leq H}\\&\times\Bigg|\sum_{\substack{r\sim R\\(r,a_2qk_1k_2)=1\\\text{cond}(\chi_2)|r}}\frac{1}{r\varphi(r)}\hat{f}\left(\frac{h}{qr}\right)e\left(\varrho_1\frac{\overline{r}}{k_1q}-\frac{a_2h}{k_1qr}\right)\Bigg|.
	\end{align*}
	Define 
	\begin{align*}
		G_1(z):=\sum_{\substack{r\leq z\\(r,a_2qk_1k_2)=1\\\text{cond}(\chi_2)|r}}\frac{r}{\varphi(r)}e\left(\varrho_1\frac{\overline{r}}{k_1q}\right),\quad
		g_2(z):=\frac{1}{z^2}\hat{f}\left(\frac{h}{qz}\right)e\left(-\frac{a_2h}{k_1qz}\right).
	\end{align*}
	For $z\in (R,2R]$, we have $g_2(z)\ll MR^{-2}$ and $g_2'(z)\ll (1+|a_2|H/(KQR))MR^{-3}\ll x^\varepsilon MR^{-3}$.
	Thus
	\begin{align*}
		&\Bigg|\sum_{\substack{r\sim R\\(r,a_2qk_1k_2)=1\\\text{cond}(\chi_2)|r}}\frac{1}{r\varphi(r)}\hat{f}\left(\frac{h}{qr}\right)e\left(\varrho_1\frac{\overline{r}}{k_1q}-\frac{a_2h}{k_1qr}\right)\Bigg|\\&\ll \max_{z\in\{R,2R\}}|G_1(z)g_2(z)|+R\max_{z\in(R,2R]}|G_1(z)g_2'(z)|\\&\ll x^\varepsilon MR^{-2}\max_{z\in[R,2R]}|G_1(z)|.
	\end{align*}
	Since $(\varrho_1,k_1q)\leq(h,k_1q)(a_2-a_{1q}k_1\overline{k_1}^{(q)},k_1q)\leq(h,k_1q)(a_2,k_1)q$,
	by Lemma 2.3, we get
	\begin{align*}
		G_1(z)\ll x^{0.1\varepsilon}\left((h,k_1q)^{1/2}(k_1q)^{1/2}+\frac{(h,k_1q)R}{k_1q}\right)(a_2,k_1)Q.
	\end{align*}
	Then
	\begin{align*}
		R_2\ll& MQ^{-1}R^{-2}x^{5.2\varepsilon}\sum_{k_1\sim K}(a_2,k_1)\sum_{k_2\sim K}\sum_{q\sim Q}\sum_{0<|h|\leq H}\Big((h,k_1q)^{1/2}(k_1q)^{1/2}+\frac{(h,k_1q)R}{k_1q}\Big)\\\ll& MQ^{-1}R^{-2}x^{5.3\varepsilon}(K^{5/2}Q^{3/2}H+HRK)\\\ll& x^{6\varepsilon}K^{5/2}Q^{3/2}R^{-1}+x^{6\varepsilon}K.\tag{3.5}
	\end{align*}
	By (3.4) and (3.5), we have
	\begin{align*}
		\mathscr{S}_2=\hat{f}(0)X_2+O(x^{6\varepsilon}K^{5/2}Q^{3/2}R^{-1}+x^{6\varepsilon}K).\tag{3.6}
	\end{align*}

	\subsection{Estimate of $\mathscr{S}_1$}Let $v=(k_1,k_2)$. Then \begin{align*}
		\mathscr{S}_1=\sum_{v\geq1}S(v)=\sum_{v\geq 1}\sum_{\substack{q\sim Q}}\sum_{\substack{r\sim R\\(r,a_2q)=1}}\sum_{(m,qr)=1}f(m)\mathop{\sum\sum}_{\substack{mk_1\equiv mk_2\equiv a_{1q}\mkern-15mu\pmod{q}\\mk_1\equiv mk_2\equiv a_2\mkern-15mu\pmod{r}\\(k_1,k_2)=v}}u_{k_1}\overline{u_{k_2}}.
	\end{align*}
	Let $\eta:=\varepsilon^2$. Suppose $Q<M$, which is not essential and we can verify it at the end. We can deduce that the contribution
	of $v> x^{\eta}$ is bounded by (We first estimate the number of $k_2$ separately. For the remaining $k_1$, we write $k_1=k$, so that $r|(mk-a_2)$.)
	\begin{align*}
		\sum_{v>x^{\eta}}S(v)\ll&\sum_{v>x^{\eta}}(Kv^{-1}Q^{-1}R^{-1}+1)x^{0.1\eta}\sum_{q\sim Q}\sum_{\substack{v|k}}|u_k|\sum_{\substack{m\sim M\\mk\equiv a_{1q}\mkern-15mu\pmod{q}\\mk\neq a_2}}\tau(mk-a_2)\\&+(KQ^{-1}R^{-1}+1)x^{0.1\eta}QR\\\ll& x^{-0.5\eta}MK^2Q^{-1}R^{-1}+x^{0.5\eta}MK+x^{0.1\eta}QR.
	\end{align*}
	Write $k_1=vd_1k_1'$ and $k_2=vk_2'$, where $d_1|v^\infty$ and $(k_1',v)=1$. Considering $v\leq x^\eta$, the contribution of $d_1> x^{\eta}$ is bounded by
	\begin{align*}
		&\ll\sum_{v\leq x^{\eta}}\sum_{\substack{x^{\eta}<d_1<2K/v\\d_1|v^\infty}}x^{0.1\eta}(Kv^{-1}Q^{-1}R^{-1}+1)\sum_{m}f(m)\sum_{q\sim Q}\sum_{\substack{k\sim K\\vd_1|k}}\sum_{\substack{m\sim M\\mk\equiv a_{1q}\mkern-15mu\pmod{q}\\mk\neq a_2}}1\\&\quad+(KQ^{-1}R^{-1}+1)x^{0.1\eta}QR\\&\ll x^{0.1\eta}\sum_{v\leq x^{\eta}}(K^2MQ^{-1}R^{-1}v^{-2}+KMv^{-1})\sum_{\substack{x^{\eta}<d_1<2K/v\\d_1|v^\infty}}\frac{1}{d_1}+(KQ^{-1}R^{-1}+1)x^{0.1\eta}QR\\&\ll x^{-0.5\eta}MK^2Q^{-1}R^{-1}+x^{0.5\eta}MK+x^{0.1\eta}QR.
	\end{align*}
	Write $\sum_{v\leq x^{\eta}}S(v)=\sum_{v\leq x^{\eta}}\sum_{d_1\leq x^\eta}S(v,d_1)+O(x^{-0.5\eta}K^2MQ^{-1}R^{-1}+x^{0.5\eta}MK+x^{0.1\eta}QR)$, where\begin{align*}
		S(v,d_1)=\sum_{\substack{q\sim Q}}\sum_{\substack{r\sim R\\(r,a_2q)=1}}\mathop{\sum\sum}_{\substack{(k_1,k_2)\in\mathcal{K}}}u_{k_1}\overline{u_{k_2}}\sum_{\substack{m\equiv a_{1q}\overline{k_1}\mkern-15mu\pmod{q}\\m\equiv a_2\overline{k_1}\mkern-15mu\pmod{r}}}f(m),
	\end{align*}
	\begin{align*}
		\mathcal{K}=\{(k_1,k_2):k_1\equiv k_2\mkern-15mu\pmod{qr},vd_1|k_1,v|k_2,(k_1,qr)=(k_1/(vd_1),v)=(k_1/v,k_2/v)=1\}.
	\end{align*} 
	By Lemma 2.6, Poisson summation yields
	\begin{align*}
		S(v,d_1)=\hat{f}(0)X_1(v,d_1)+R(v,d_1)+O(1),
	\end{align*}
	where
	\begin{align*}
		&X_1(v,d_1)=\sum_{\substack{q\sim Q}}\sum_{\substack{r\sim R\\(r,a_2q)=1}}\frac{1}{qr}\mathop{\sum\sum}_{(k_1,k_2)\in\mathcal{K}}u_{k_1}\overline{u_{k_2}},\\&R(v,d_1)=\sum_{\substack{q\sim Q}}\sum_{\substack{r\sim R\\(r,a_2q)=1}}\frac{1}{qr}\mathop{\sum\sum}_{(k_1,k_2)\in\mathcal{K}}u_{k_1}\overline{u_{k_2}}\sum_{0<|h|\leq H}\hat{f}\left(\frac{h}{qr}\right)e\left(-a_{1q}h\frac{\overline{k_1r}}{q}-a_2h\frac{\overline{k_1q}}{r}\right).
	\end{align*}
	Following a similar argument, we have
	\begin{align*}
		\mathop{\sum\sum}_{v>x^{\eta}\ \text{or}\ d_1>x^{\eta}}\hat{f}(0)X_1(v,d_1)\ll x^{-0.5\eta}MK^2Q^{-1}R^{-1}+x^{0.5\eta}MK+x^{0.1\eta}QR.
	\end{align*}
	Thus we arrive at
	\begin{align*}
		\mathscr{S}_1=\hat{f}(0)X_1+O\Bigg(\sum_{v\leq x^\eta}\sum_{d_1\leq x^{\eta}}|R(v,d_1)|\Bigg)+T_1,\tag{3.7}
	\end{align*}
	where
	\begin{align*}
		X_1&=\sum_{\substack{q\sim Q}}\sum_{\substack{r\sim R\\(r,a_2q)=1}}\frac{1}{qr}\mathop{\sum\sum}_{\substack{k_1\equiv k_2\mkern-15mu\pmod{qr}\\(k_1k_2,qr)=1}}u_{k_1}\overline{u_{k_2}},\\T_1&\ll x^{-0.5\varepsilon^2} MK^2Q^{-1}R^{-1}+x^{0.5\varepsilon}MK+x^{0.1\varepsilon}QR.
	\end{align*}
	
	The estimation of $R(v,d_1)$ is the most important.
	We claim that
	\begin{proposition}	
		Given $d_1|v^\infty$, $d_2|v$,  $|vd_1d_2|\leq x^{\varepsilon}$, $|a_2|\leq x^{1+\varepsilon}$, let $Q$, $R$, $K$, $M$, $N$, $L$, $H>1$ be real numbers satisfying\begin{align*}
			x^{1-\varepsilon}<MNL\leq x,\quad H\ll x^{\varepsilon} QRM^{-1},\quad d_1K\asymp d_2NL,
		\end{align*}
		\begin{align*}
			x^{-20\varepsilon}NL<QR\ll x^{-\varepsilon}d_2NL,\quad  Q\ll x^{\varepsilon}N.
		\end{align*} For all $k_1\sim K$, $n\sim N$, $l\sim L$, all sequences $(u_{k_1})$, $(\beta_n)$, $(\gamma_l)$ with
		\begin{align*}
			|u_{k_1}|\leq \tau(k_1),\quad |\beta_n|\leq 1,\quad |\gamma_l|\leq 1,
		\end{align*}
		and $u_{k_1}=0$ whenever $(k_1,v)>1$, and for all smooth functions $f$ supported on $[M/2,3M]$ satisfying $f(m)\geq 0$ and $f^{(j)}(m)\ll_j M^{-j}$ for $j\geq 0$, define
		\begin{align*}
			\mathcal{R}:=&\sum_{\substack{q\sim Q}}\sum_{\substack{r\sim R\\(r,a_2q)=1}}\frac{1}{qr}\mathop{\sum\sum}_{\substack{d_1k_1\equiv d_2nl\mkern-15mu\pmod{qr}\\(d_1k_1,d_2nl)=1}}u_{k_1}\beta_n\gamma_l\sum_{0<|h|\leq H}\hat{f}\left(\frac{h}{qr}\right)\\&\times e\left(-a_{1q}h\frac{\overline{vd_1k_1r}}{q}-a_2h\frac{\overline{vd_1k_1q}}{r}\right).
		\end{align*}
		Suppose\begin{align*}
			Q^{1/2}N^{1/2}L<x^{1/2-O(\varepsilon)},\ QN^{2}L^{7/4}<x^{1-O(\varepsilon)},
		\end{align*}
		then one has
		\begin{align*}
			MQR\cdot \mathcal{R}\ll x^{2-\varepsilon}.
		\end{align*}
	\end{proposition}
	
	\begin{remark}
		In Proposition 3.1, we suppress the condition $(v,qr)=1$ for brevity. We also have $(d_1d_2k_1nl,qr)=1$, since $(d_1k_1,d_2nl)=1$. In Proposition 3.1, we allow $N$ and $L$ to differ from the original $N$ and $L$ by a factor of $x^\varepsilon$, due to $v$, $d_1$ and $d_2$. But this has no effect. 
	\end{remark}

	\subsubsection{Proof of Proposition 3.1}
	
	For convenience, we recall that
	\begin{align*}
		A\pp B
	\end{align*}
	means $A\ll_\varepsilon x^{O(\varepsilon)}B$.
	Let $\nu:=vd_1$. Without ambiguity, we write\begin{align*}
		u_{k_2}=\sum_{nl=k_2}\beta_n\gamma_l.
	\end{align*}
	Then
	\begin{align*}
		\mathcal{R}=\sum_{\substack{q\sim Q}}\sum_{\substack{r\sim R\\(r,a_2q)=1}}\frac{1}{qr}\mathop{\sum\sum}_{\substack{d_1k_1\equiv d_2k_2\mkern-15mu\pmod{qr}\\(d_1k_1,d_2k_2)=1}}u_{k_1}u_{k_2}\sum_{0<|h|\leq H}\hat{f}\left(\frac{h}{qr}\right)e\left(-a_{1q}h\frac{\overline{\nu k_1r}}{q}-a_2h\frac{\overline{\nu k_1q}}{r}\right).
	\end{align*}
	Since
	\begin{align*}
		\nu k_1q\overline{\nu k_1q}^{(r)}+\nu(d_2k_2-d_1k_1)\overline{\nu d_2k_2}^{(k_1)}+k_1qr\overline{k_1qr}^{(\nu)}+\nu k_1r\overline{\nu k_1r}^{(  q)}\equiv1\mkern-15mu\pmod{\nu k_1qr},
	\end{align*}
	and \begin{align*}
		\frac{\overline{\nu d_2k_2}}{k_1}+\frac{\overline{k_1}}{\nu d_2k_2}\equiv\frac{1}{\nu d_2k_1k_2}\mkern-15mu\pmod{1},
	\end{align*}
	we have
	\begin{align*}
		-\frac{\overline{\nu k_1q}}{r}&\equiv\frac{d_2k_2-d_1k_1}{qr}\frac{\overline{\nu d_2k_2}}{k_1}+\frac{\overline{k_1qr}}{\nu}+\frac{\overline{\nu k_1r}}{q}-\frac{1}{\nu k_1 qr}\\&\equiv-\frac{d_2k_2-d_1k_1}{qr}\frac{\overline{k_1}}{\nu d_2k_2}+\frac{\overline{k_1qr}}{\nu}+\frac{\overline{\nu k_1r}}{q}+\frac{d_2k_2-d_1k_1}{\nu d_2k_1k_2qr}-\frac{1}{\nu k_1 qr} \mkern-15mu\pmod{1}.
	\end{align*}
	Then \begin{align*}
		e\left(-a_{1q}h\frac{\overline{\nu k_1r}}{q}-a_2h\frac{\overline{\nu k_1q}}{r}\right)=&e\Bigg(-a_2h\frac{d_2k_2-d_1k_1}{qr}\frac{\overline{k_1}}{\nu d_2k_2}+a_2h\frac{\overline{k_1qr}}{\nu}+(a_2-a_{1q})h\frac{\overline{\nu k_1r}}{q}\\&+a_2h\frac{d_2k_2-d_1k_1}{\nu d_2k_1k_2qr}-\frac{a_2h}{\nu k_1 qr}\Bigg).
	\end{align*}
	Let $t:=(d_2k_2-d_1k_1)/(qr)$. Then we have $0\neq t\pp1$. Without loss of generality, we assume $t>0$. 
	Then we get
	\begin{align*}
		\mathcal{R}\pp&\max_{0\neq t\ll x^{O(\varepsilon)}}\Bigg|\sum_{\substack{q\sim Q}}\sum_{\substack{k_1,k_2\\(d_1k_1,d_2k_2)=1\\d_2k_2\equiv d_1k_1\mkern-15mu\pmod{qt}\\r\sim R\\(r,a_2)=1}}u_{k_1}u_{k_2}\sum_{0<|h|\leq H}\frac{t}{d_2k_2-d_1k_1}\hat{f}\left(\frac{ht}{d_2k_2-d_1k_1}\right)\\&\times e\left(-a_2ht\frac{\overline{k_1}}{\nu d_2k_2}+\frac{a_2ht}{\nu d_2k_1k_2}-\frac{a_2h}{\nu k_1 qr}+a_2h\frac{\overline{k_1qr}}{\nu}+(a_2-a_{1q})h\frac{\overline{\nu k_1r}}{q}\right)\Bigg|,
	\end{align*}
	where $r=(d_2k_2-d_1k_1)/(qt)$.
	By means of the Fourier transform, one can obtain\begin{align*}
		\frac{t}{d_2k_2-d_1k_1}\hat{f}\left(\frac{ht}{d_2k_2-d_1k_1}\right)=\int_{0}^{3MQ^{-1}R^{-1}}\int_{-\infty}^{\infty}t\hat{f}(\eta t)e(-\eta\xi(d_2k_2-d_1k_1))e(\xi h)\mathrm{d}\eta\mathrm{d}\xi.
	\end{align*}
	Then $\mathcal{R}$ becomes
	\begin{align*}
		\mathcal{R}\pp&\max_{0\neq t\ll x^{O(\varepsilon)}} \sum_{\substack{q\sim Q}}\sum_{k_1\sim K}\sum_{l\sim L }\int_{0}^{3MQ^{-1}R^{-1}}\int_{-\infty}^{\infty}
		|t\hat{f}(\eta t)|\\&\times\Bigg|\sum_{0<|h|\leq H}\sum_{\substack{n\\(d_2nl,d_1  k_1)=1\\d_2nl\equiv d_1k_1\mkern-15mu\pmod{qt}\\r\sim R\\(r,a_2)=1}}\beta_ne\Bigg(-a_2ht\frac{\overline{k_1}}{\nu d_2nl}+\frac{a_2ht}{\nu d_2k_1nl}-\frac{a_2h}{\nu k_1 qr}
		\\& +a_2h\frac{\overline{k_1qr}}{\nu}+(a_2-a_{1q})h\frac{\overline{\nu k_1r}}{q}-d_2nl\eta\xi+h\xi \Bigg)\Bigg|\mathrm{d}\eta\mathrm{d}\xi.\tag{3.8}
	\end{align*}
	We perform the change of variables $\eta/l\leftarrow\eta$,  which yields
	\begin{align*}
		\mathcal{R}\pp& \max_{0\neq t\ll x^{O(\varepsilon)}}\sum_{\substack{q\sim Q}} \sum_{k_1\sim K}\sum_{l\sim L}\int_{0}^{3MQ^{-1}R^{-1}}\int_{-\infty}^{\infty}\Big|\frac{t}{l}\hat{f}\left(\frac{\eta t}{l}\right)\Big|\\&\times\Bigg|\sum_{0<|h|\leq H} \sum_{\substack{n\\(d_2nl,d_1 k_1)=1\\d_2nl\equiv d_1k_1\mkern-15mu\pmod{qt}\\r\sim R\\(r,a_2)=1}}\beta_ne\Bigg(-a_2ht\frac{\overline{k_1}}{\nu d_2nl}+\frac{a_2ht}{\nu d_2k_1nl}-\frac{a_2h}{\nu k_1 qr}\\&+a_2h\frac{\overline{k_1qr}}{\nu}+(a_2-a_{1q})h\frac{\overline{\nu k_1r}}{q}-d_2n\eta\xi+h\xi \Bigg)\Bigg|\mathrm{d}\eta\mathrm{d}\xi ,
	\end{align*}
	where
	\begin{align*}
		\Big|\frac{t}{l}\hat{f}\left(\frac{\eta t}{l}\right)\Big|\ll\min\{tML^{-1},|\eta|^{-1},LM|\eta|^{-2}\}.
	\end{align*}
	Then we arrive
	\begin{align*}
		\mathcal{R}\pp& MQ^{-1}R^{-1}\sup_{\substack{\eta\in(-\infty,\infty)\\\xi\in[0,3MQ^{-1}R^{-1}]}}\max_{0\neq t\ll x^{O(\varepsilon)}} \sum_{\substack{q\sim Q}}\sum_{k_1\sim K}\sum_{l\sim L}\Bigg|\sum_{0<|h|\leq H} \sum_{\substack{n\\(d_2nl,d_1 k_1)=1\\d_2nl\equiv d_1k_1\mkern-15mu\pmod{qt}\\r\sim R\\(r,a_2)=1}}\beta_n\\&\times e\Bigg(-a_2ht\frac{\overline{k_1}}{\nu d_2nl}+\frac{a_2ht}{\nu d_2k_1nl}-\frac{a_2h}{\nu k_1 qr}\\&+a_2h\frac{\overline{k_1qr}}{\nu}+(a_2-a_{1q})h\frac{\overline{\nu k_1r}}{q}-d_2n\eta\xi+h\xi \Bigg)\Bigg|.
	\end{align*}
	By Cauchy-Schwarz, we have
	\begin{align*}
		\mathcal{R}\pp MQ^{-1}R^{-1}\sup_{\substack{\eta\in(-\infty,\infty)\\\xi\in[0,3MQ^{-1}R^{-1}]}}\max_{0\neq t\ll x^{O(\varepsilon)}}(QKL)^{1/2}\mathcal{R}'^{1/2},\tag{3.9}
	\end{align*}
	where
	\begin{align*}
		\mathcal{R}'=&\sum_{\substack{q\sim Q}}\sum_{k_1\sim K}\sum_{l\sim L} \Bigg|\sum_{0<|h|\leq H}\sum_{\substack{n\\(d_2nl,d_1 k_1)=1\\d_2nl\equiv d_1k_1\mkern-15mu\pmod{qt}\\r\sim R\\(r,a_2)=1}}\beta_ne\Bigg(-a_2ht\frac{\overline{k_1}}{\nu d_2nl}+\frac{a_2ht}{\nu d_2k_1nl}-\frac{a_2h}{\nu k_1 qr}\\& +a_2h\frac{\overline{k_1qr}}{\nu}+(a_2-a_{1q})h\frac{\overline{\nu k_1r}}{q}-d_2n\eta\xi+h\xi \Bigg)\Bigg|^2\\\leq& \sum_{\substack{q\sim Q}}\sum_{l\sim L}\sum_{0<|h_1|\leq H}\sum_{0<|h_2|\leq H}\sum_{n_1\sim N}\sum_{n_2\sim N}  \Bigg|\sum_{\substack{k_1\sim K\\(d_2n_1n_2l,d_1 k_1)=1\\d_1k_1\equiv d_2n_1l\equiv d_2n_2l\mkern-15mu\pmod{qt}\\r_1,r_2\sim R\\(r_1r_2,a_2)=1}}\\&\times e\Bigg(a_2t\frac{\overline{k_1}}{\nu d_2l[n_1,n_2]}\frac{n_1h_2-n_2h_1}{(n_1,n_2)}+a_2t\frac{n_2h_1-n_1h_2}{\nu d_2k_1n_1n_2l}-\frac{a_2h_1}{\nu k_1qr_1}+\frac{a_2h_2}{\nu k_1qr_2}\\&+a_2h_1\frac{\overline{k_1qr_1}}{\nu}-a_2h_2\frac{\overline{k_1qr_2}}{\nu}+(a_2-a_{1q})h_1\frac{\overline{\nu k_1r_1}}{q}-(a_2-a_{1q})h_2\frac{\overline{\nu k_1r_2}}{q} \Bigg)\Bigg|,
	\end{align*}
	with $r_1=(d_2n_1l-d_1k_1)/(qt)$ and $r_2=(d_2n_2l-d_1k_1)/(qt)$.
	The contribution of $n_1h_2=n_2h_1$ for $R'$ is 
	\begin{align*}
		\mathcal{R}'(n_1h_2=n_2h_1)\pp HKNL.\tag{3.10}
	\end{align*} 
	For $n_1h_2\neq n_2h_1$, we rewrite $(r_1r_2,a_2)=1$ as
	\begin{align*}
		\sum_{\sigma_1|a_2}\mu(\sigma_1)1_{r_1\equiv 0\mkern-15mu\pmod{\sigma_1}}\sum_{\sigma_2|a_2}\mu(\sigma_2)1_{r_2\equiv 0\mkern-15mu\pmod{\sigma_2}}
	\end{align*} and consider
	\begin{align*}
		&r_1\equiv w_1\mkern-15mu\pmod{\nu },\\&r_1\equiv w_2\mkern-15mu\pmod{q},
	\end{align*}where $0<w_1<\nu $, $(w_1,\nu )=1$ and $0<w_2<q$, $(w_2,q)=1$. In fact, we can rewrite the congruence conditions as
	\begin{align*}
		&d_2n_1l-d_1k_1\equiv 0\mkern-15mu\pmod{\sigma_1 qt},\\& d_2n_2l-d_1k_1\equiv 0\mkern-15mu\pmod{\sigma_2 qt},\\
		&d_1k_1\equiv d_2n_1l-w_1qt\mkern-15mu\pmod{\nu qt},\\&d_1k_1\equiv d_2n_1l-w_2qt\mkern-15mu\pmod{q^2t}.
	\end{align*}
	We have now fixed $k_1$ modulo $q$. We also have\begin{align*}
		&qtr_2=d_2n_2l-d_1k_1\equiv d_2n_2l-d_2n_1l+w_1qt\mkern-15mu\pmod{\nu  qt},\\&qtr_2\equiv d_2n_2l-d_2n_1l+w_2qt\mkern-15mu\pmod{q^2t},
	\end{align*}
	which implies that
	\begin{align*}
		&n_2\equiv n_1\mkern-15mu\pmod{qt},\\&r_2\equiv (d_2n_2l-d_2n_1l)/(qt)+w_1\mkern-15mu\pmod{\nu },\\&r_2\equiv (d_2n_2l-d_2n_1l)/(qt)+w_2\mkern-15mu\pmod{q}.
	\end{align*}
	Thus we have
	\begin{align*}
		\mathcal{R}'&(n_1h_2\neq n_2h_1)\\\ll&\sum_{\substack{q\sim Q}}\sum_{l\sim L}\sum_{0<|h_1|\leq H}\sum_{0<|h_2|\leq H}\mathop{\sum_{n_1\sim N}\sum_{n_2\sim N}}_{\substack{n_1\equiv n_2\mkern-15mu\pmod{qt}\\n_1h_2\neq n_2h_1}}\sum_{\substack{0<w_1<\nu \\(w_1,\nu )=1}}\sum_{\substack{0<w_2<q\\(w_2,q)=1}}\sum_{\sigma_1|a_2}\sum_{\sigma_2|a_2}\\&\times\Bigg|  \sum_{\substack{k_1\sim K\\r_1,r_2\sim R\\(d_2n_1n_2l,d_1k_1)=1\\d_1k_1\equiv d_2n_1l\mkern-15mu\pmod{\sigma_1 qt}\\ d_1k_1\equiv d_2n_2l\mkern-15mu\pmod{\sigma_2 qt}\\d_1k_1\equiv d_2n_1l-w_1qt\mkern-15mu\pmod{\nu qt}\\d_1k_1\equiv d_2n_1l-w_2qt\mkern-15mu\pmod{q^2t}
		}} e\Bigg(a_2t\frac{\overline{k_1}}{\nu d_2l[n_1,n_2]}\frac{n_1h_2-n_2h_1}{(n_1,n_2)}\\&+a_2t\frac{n_2h_1-n_1h_2}{\nu d_2k_1n_1n_2l}-\frac{a_2h_1}{\nu k_1qr_1}+\frac{a_2h_2}{\nu k_1qr_2} \Bigg)\Bigg|.\tag{3.11}
	\end{align*}
	For $z\sim K$ and $R<(d_2n_1l-d_1z)/(qt), (d_2n_2l-d_1z)/(qt)\leq 2R$, let \begin{align*}
		g_1(z):=e\left(a_2t\frac{n_2h_1-n_1h_2}{\nu d_2zn_1n_2l}-\frac{a_2h_1}{\nu zqr_1}+\frac{a_2h_2}{\nu z qr_2}\right),
	\end{align*}
	where $r_1=(d_2n_1l-d_1z)/(qt)$ and $r_2=(d_2n_2l-d_1z)/(qt)$. Then\begin{align*}
		g_1'(z)\pp |a_2|\left(\frac{H}{NLK^2}+\frac{H}{QRK^2}+\frac{H}{Q^2R^2K}\right)\pp K^{-1}.\tag{3.12}
	\end{align*} 
	Then using (3.12) and partial summation, we get
	\begin{align*}
		&\Bigg|\sum_{\substack{k_1\sim K\\r_1,r_2\sim R\\(d_2n_1n_2l,d_1k_1)=1\\d_1k_1\equiv d_2n_1l\mkern-15mu\pmod{\sigma_1 qt}\\ d_1k_1\equiv d_2n_2l\mkern-15mu\pmod{\sigma_2 qt}\\d_1k_1\equiv d_2n_1l-w_1qt\mkern-15mu\pmod{\nu qt}\\d_1k_1\equiv d_2n_1l-w_2qt\mkern-15mu\pmod{q^2t}
		}}e\left(\frac{a_2t\overline{k_1}}{\nu d_2l[n_1,n_2]}\frac{n_1h_2-n_2h_1}{(n_1,n_2)}\right)g(k_1)\Bigg|\\&\pp\max_{K<K'\leq 2K}\Bigg|\sum_{\substack{K<k_1\leq K'\\r_1,r_2\sim R\\(d_2n_1n_2l,d_1k_1)=1\\d_1k_1\equiv d_2n_1l\mkern-15mu\pmod{\sigma_1 qt}\\ d_1k_1\equiv d_2n_2l\mkern-15mu\pmod{\sigma_2 qt}\\d_1k_1\equiv d_2n_1l-w_1qt\mkern-15mu\pmod{\nu qt}\\d_1k_1\equiv d_2n_1l-w_2qt\mkern-15mu\pmod{q^2t}
		}}e\left(\frac{a_2t\overline{k_1}}{\nu d_2l[n_1,n_2]}\frac{n_1h_2-n_2h_1}{(n_1,n_2)}\right)\Bigg|
	\end{align*}
	Since $K\ll \nu d_2l[n_1,n_2]$, by Lemma 2.3, we get\begin{align*}
		&\Bigg|\sum_{\substack{K<k_1\leq K'\\r_1,r_2\sim R\\(d_2n_1n_2l,d_1k_1)=1\\d_1k_1\equiv d_2n_1l\mkern-15mu\pmod{\sigma_1 qt}\\ d_1k_1\equiv d_2n_2l\mkern-15mu\pmod{\sigma_2 qt}\\d_1k_1\equiv d_2n_1l-w_1qt\mkern-15mu\pmod{\nu qt}\\d_1k_1\equiv d_2n_1l-w_2qt\mkern-15mu\pmod{q^2t}
		}}e\left(\frac{a_2t\overline{k_1}}{\nu d_2l[n_1,n_2]}\frac{n_1h_2-n_2h_1}{(n_1,n_2)}\right)\Bigg|\\&\pp  (a_2t(n_1h_2-n_2h_1)/(n_1,n_2),\nu d_2l[n_1,n_2])^{1/2}(\nu d_2l[n_1,n_2])^{1/2}\\&\pp (n_1n_2l)^{1/2}(a_2(n_1h_2-n_2h_1)/(n_1,n_2),l[n_1,n_2]/(n_1,n_2))^{1/2}\\&\pp (n_1n_2l)^{1/2}(a_2(n_1h_2-n_2h_1)/(n_1,n_2),l)(a_2h_1,n_1)^{1/2}(a_2h_2,n_2)^{1/2}.
	\end{align*}
	Summing over all variables, we have
	\begin{align*}
		\mathcal{R}'(n_1h_2\neq n_2h_1)\pp&NL^{1/2}\sum_{q\sim Q}\sum_{\substack{0<w_1<\nu \\(w_1,\nu )=1}}\sum_{\substack{0<w_2<q\\(w_2,q)=1}}\sum_{n_1\sim N}\sum_{\substack{n_2\sim N\\n_1\equiv n_2\mkern-15mu\pmod{ qt}}}\sum_{0<|h_1|\leq H}\sum_{0<|h_2|\leq H}\\&\times(a_2h_1,n_1)^{1/2}(a_2h_2,n_2)^{1/2}\sum_{l\sim L}(a_2(n_1h_2-n_2h_1)/(n_1,n_2),l)\\\pp& QH^2NL^{3/2}\sum_{n_1\sim N}\sum_{\substack{n_2\sim N}}(a_2,n_1)^{1/2}(a_2,n_2)^{1/2}\sum_{\substack{q|(n_1-n_2)\\q\sim Q}}1\\\pp&QH^2N^3L^{3/2}+Q^2H^2N^2L^{3/2}\\\pp&QH^2N^3L^{3/2}.\tag{3.13}
	\end{align*}
	We conclude from (3.9), (3.10), (3.13) that
	\begin{align*}
		MQR\cdot\mathcal{R}&\pp MQRMQ^{-1}R^{-1}(QKL)^{1/2}R'^{1/2}\\&\pp
		M^2(QKL)^{1/2}(HKNL+QH^2N^3L^{3/2})^{1/2}\\&\pp
		Q^{1/2}H^{1/2}KM^2N^{1/2}L+QHK^{1/2}M^2N^{3/2}L^{5/4}	\\&\pp Q^{1/2}M^{3/2}N^{2}L^{5/2}+QMN^{3}L^{11/4}.
	\end{align*}
	Suppose\begin{align*}
		Q^{1/2}N^{1/2}L<x^{1/2-O(\varepsilon)},\ QN^{2}L^{7/4}<x^{1-O(\varepsilon)},
	\end{align*}
	then we have
	\begin{align*}
		MQR\cdot\mathcal{R}\ll x^{2-\varepsilon}.
	\end{align*}
	We complete the proof of Proposition 3.1.
	\qed
	
	Now we assume that all conditions in Proposition 3.1 are satisfied, and then estimate $\mathscr{S}_1$.
	Consider $\max(v,d_1)\leq x^\eta$ with $\eta=\varepsilon^2$. Write $k_1=vd_1k_1''$ and $k_2=vd_2k_2''$, where $d_2=(v,k_2/v)$. By $u_{k_2}=\sum_{nl=vd_2k_2''}\beta_n\gamma_l$, define
	\begin{align*}
		\delta_2:=(l,v),\quad \delta_1:=v/\delta_2,\quad l':=l/\delta_2,\quad n':=n/\delta_1. 
	\end{align*}
	Then we have
	\begin{align*}
		u_{k_2}=\sum_{nl=vd_2k_2''}\beta_n\gamma_l=\sum_{\delta_1\delta_2=v}\sum_{\substack{n'l'=d_2k_2''\\(l',\delta_1)=1}}\beta_{\delta_1n'}\gamma_{\delta_2l'}=\sum_{\delta_1\delta_2=v}\sum_{\substack{\delta_3\delta_4=d_2\\(\delta_1,\delta_4)=1}}\sum_{\substack{n''l''=k_2''\\(\delta_1\delta_3,l'')=1}}\beta_{\delta_1\delta_3n''}\gamma_{\delta_2\delta_4l''}.
	\end{align*}
	Write $R(v,d_1)=\sum_{d_2|v,(d_2,d_1)=1}R(v,d_1,d_2)$, where $R(v,d_1,d_2)$ is defined by
	\begin{align*}
		R(v,d_1,d_2):=\sum_{\delta_1\delta_2=v}\sum_{\substack{\delta_3\delta_4=d_2\\(\delta_1,\delta_4)=1}}R(v,d_1,d_2,\delta_1,\delta_2,\delta_3,\delta_4),
	\end{align*}
	with
	\begin{align*}
		R(v,d_1,d_2,\delta_1,\delta_2,\delta_3,\delta_4)=&\sum_{\substack{q\sim Q\\(q,v)=1}}\sum_{\substack{r\sim R\\(r,va_2q)=1}}\frac{1}{qr}\mathop{\sum\sum}_{\substack{d_1k_1''\equiv d_2n''l''\mkern-15mu\pmod{qr}\\(d_1k_1'',d_2n''l'')=1}}u_{k_1''}'\beta_{n''}'\gamma_{l''}'\\&\times\sum_{0<|h|\leq H}\hat{f}\left(\frac{h}{qr}\right) e\left(-a_{1q}h\frac{\overline{vd_1k_1''r}}{q}-a_2h\frac{\overline{vd_1k_1''q}}{r}\right),
	\end{align*}
	\begin{align*}
		u'_{k_1''}=1_{(k_1'',v)=1}u_{vd_1k_1''},\ \beta_{n''}'=1_{(n'',v/d_2)=1}\overline{\beta_{\delta_1\delta_3n''}},\ \gamma_{l''}'=1_{(l'',v/d_2)=1}\overline{\gamma_{\delta_2\delta_4l''}}.
	\end{align*}
	By Proposition 3.1, we have
	\begin{align*}
		R(v,d_1,d_2,\delta_1, \delta_2,\delta_3,\delta_4)\ll x^{2-\varepsilon}M^{-1}Q^{-1}R^{-1}.
	\end{align*}
	Summing over all variables, we get
	\begin{align*}
		\mathscr{S}_1=&\hat{f}(0)X_1+O(x^{2-\varepsilon/2}M^{-1}Q^{-1}R^{-1}+x^{-0.5\varepsilon^2}MK^2Q^{-1}R^{-1}+x^{0.5\varepsilon}MK+x^{0.1\varepsilon}QR).\tag{3.14}
	\end{align*}

	\subsection{Estimation of the main term}
	Changing the order of summation, we have
	\begin{align*}
		X_1&-2Re X_2+X_3=\sum_{\substack{q\sim Q}}\sum_{\substack{r\sim R\\(r,a_2q)=1}}\frac{1}{qr}\sum_{\substack{0<b<qr\\(b,qr)=1}}\Bigg|\frac{1}{\varphi(qr)}\sum_{(k,qr)=1}u_k\omega_\varepsilon(kb,q)\omega_\varepsilon(kb,r)\Bigg|^2\\&-2Re \sum_{\substack{q\sim Q}}\sum_{\substack{r\sim R\\(r,a_2q)=1}}\frac{1}{qr\varphi(qr)}\mathop{\sum_{k_1}\sum_{k_2}}_{(k_1k_2,qr)=1}u_{k_1}\overline{u_{k_2}}\omega_\varepsilon(k_1\overline{k_2},q)\omega_\varepsilon(k_1\overline{k_2},r)\\&+\sum_{\substack{q\sim Q}}\sum_{\substack{r\sim R\\(r,a_2q)=1}}\frac{1}{qr}\mathop{\sum\sum}_{\substack{k_1\equiv k_2\mkern-15mu\pmod{qr}\\(k_1k_2,qr)=1}}u_{k_1}\overline{u_{k_2}}\\=&\sum_{\substack{q\sim Q}}\sum_{\substack{r\sim R\\(r,a_2q)=1}}\frac{1}{qr}\sum_{\substack{0<b<qr\\(b,qr)=1}}\Bigg|\sum_{k\equiv b\mkern-15mu\pmod{qr}}u_k-\frac{1}{\varphi(qr)}\sum_{(k,qr)=1}u_k\omega_\varepsilon(k\overline{b},q)\omega_\varepsilon(k\overline{b},r)\Bigg|^2\\\leq&\sum_{\substack{q\sim Q}}\sum_{\substack{r\sim R\\(r,a_2q)=1}}\frac{1}{qr\varphi(qr)}\sum_{\substack{\chi\ \text{primitive}\\\text{cond}(\chi)>x^{\varepsilon}\\\text{cond}(\chi)|qr}}\Bigg|\sum_{(k,qr)=1}u_k\chi(k)\Bigg|^2\\=&\frac{x^{0.1\varepsilon}}{QR}\sum_{QR<s\leq 4QR}\frac{1}{\varphi(s)}\sum_{\substack{\chi\ \text{primitive}\\\text{cond}(\chi)>x^{\varepsilon}\\\text{cond}(\chi)|s}}\Bigg|\sum_{(k,s)=1}u_k\chi(k)\Bigg|^2.
	\end{align*}
	By partial summation, we have 
	\begin{align*}
		X_1-2Re X_2+X_3&\ll\frac{x^{0.1\varepsilon}}{QR}\sum_{QR<s\leq 4QR}\frac{1}{\varphi(s)}\sum_{\substack{\chi\ \text{primitive}\\\text{cond}(\chi)>x^{\varepsilon}\\\text{cond}(\chi)|s}}\Bigg|\sum_{(k,s)=1}u_k\chi(k)\Bigg|^2\\&\ll\frac{x^{0.1\varepsilon}}{QR}\sum_{x^{\varepsilon}<d\leq 4QR}\sum_{\substack{\chi\mkern-15mu\pmod{d}\\\chi\  \text{primitive}}}\sum_{\substack{s\leq 4QR\\d|s}}\frac{1}{\varphi(s)}\Bigg|\sum_{(k,s/d)=1}u_k\chi(k)\Bigg|^2
		\\&\ll\frac{x^{0.1\varepsilon}}{QR}\sum_{x^{\varepsilon}<d\leq 4QR}\sum_{\substack{\chi\mkern-15mu\pmod{d}\\\chi\  \text{primitive}}}\sum_{\substack{s\leq 4QR\\d|s}}\frac{\tau(s/d)}{\varphi(s)}\sum_{e|(s/d)}\Bigg|\sum_{e|k}u_k\chi(k)\Bigg|^2\\&\ll\frac{x^{0.2\varepsilon}}{QR}\sum_{e\leq 4QR}\frac{1}{\varphi(e)}\int_{x^\varepsilon}^{\infty}\sum_{x^{\varepsilon}<d\leq \min\{4QR,t\}}\frac{d}{\varphi(d)}\sum_{\substack{\chi\mkern-15mu\pmod{d}\\\chi\  \text{primitive}}}\Bigg|\sum_{k}u_{ek}\chi(k)\Bigg|^2\frac{\mathrm{d}t}{t^2}
	\end{align*}
	We apply Lemma 2.4 to get
	\begin{align*}
		X_1-2Re X_2+X_3&\ll\frac{x^{0.3\varepsilon}}{QR}\sum_{e\leq 4QR}\frac{1}{\varphi(e)}\left(x^{-\varepsilon}\left(\frac{K}{e}+1\right)+QR\right)\left(\frac{K}{e}+1\right)\\&\ll x^{-0.5\varepsilon}K^2Q^{-1}R^{-1}+x^{0.5\varepsilon}K.
	\end{align*}
	Suppose $QR\ll x^{-\varepsilon}NL$. By $MNL=x$, the contribution of $X_1-2Re X_2+X_3$ for $\Delta^2$ is\begin{align*}
		&\ll MQRM(x^{-0.5\varepsilon}K^2Q^{-1}R^{-1}+x^{0.5\varepsilon}K)\\&\ll x^{-0.5\varepsilon}M^2N^2L^2+x^{0.5\varepsilon}QRM^2NL\\&\ll x^{2-0.5\varepsilon}.\tag{3.15}
	\end{align*}

	\subsection{Conclusion}
	With the conditions of Proposition 3.1, we conclude from (3.2), (3.3), (3.6), (3.14), (3.15) that
	\begin{align*}
		&\Delta^2\\&\ll MQR(x^{12\varepsilon}K^2Q^{-1}R^{-1}+x^{6\varepsilon}K^{5/2}Q^{3/2}R^{-1}+x^{6\varepsilon}K
		\\&\quad+x^{2-0.5\varepsilon}M^{-1}Q^{-1}R^{-1}+x^{-0.5\varepsilon^2} MK^2Q^{-1}R^{-1}+x^{0.5\varepsilon}MK+x^{0.1\varepsilon}QR)\\&\ll
		MQR(x^{12\varepsilon}K^2Q^{-1}R^{-1}+x^{6\varepsilon}K^{5/2}Q^{3/2}R^{-1}+x^{-0.5\varepsilon^2} MK^2Q^{-1}R^{-1}+x^{6\varepsilon}MK
		)+x^{2-0.5\varepsilon^2}\\&\ll x^{12\varepsilon}MN^2L^2+ x^{6\varepsilon}Q^{5/2}MN^{5/2}L^{5/2}+x^{6\varepsilon}QRM^2NL+x^{2-0.5\varepsilon^2}.
	\end{align*}
	Suppose
	\begin{align*}
		Q^{1/2}N^{1/2}L<x^{1/2-100\varepsilon},\ QN^{2}L^{7/4}<x^{1-100\varepsilon},
	\end{align*}
	\begin{align*}
		Q^{5/2}N^{3/2}L^{3/2}<x^{1-7\varepsilon},\ x^{-13\varepsilon}NL<QR<x^{-7\varepsilon}NL,\ Q\ll \min(x^\varepsilon N,M).
	\end{align*}
	then we have
	\begin{align*}
		\Delta^2\ll x^{2-0.5\varepsilon^2}.
	\end{align*}
	In fact, by $QR<x^{-7\varepsilon}NL$, we have
	\begin{align*}
		&QR^{1/2}L^{1/2}=(QR)^{1/2}Q^{1/2}L^{1/2}<x^{-3.5\varepsilon}Q^{1/2}N^{1/2}L<x^{1/2-103.5\varepsilon},\\&Q^3R^2L^{-1/4}=(QR)^2QL^{-1/4}< x^{-14\varepsilon}QN^2L^{7/4}<x^{1-114\varepsilon}.
	\end{align*}
	Thus the best result for $R$ is 
	\begin{align*}
		L\approx x^{4/9}Q^{-4/9},\  R=x^{5/9-310\varepsilon}Q^{-14/9}.
	\end{align*}
	Let
	\begin{align*}
		M_0:=x^{4/9+300\varepsilon}Q^{5/9},\ L_0:=x^{4/9-300\varepsilon}Q^{-4/9},\ N_0:=x^{1/9}Q^{-1/9},
	\end{align*}
	For $y$ as in Theorem 1.1, by Lemma 2.7, any resulting values of $M$, $N$, $L$ with\begin{align*}
		M_0\leq M\leq y\frac{M_0}{2},\ L_0\leq L\leq y\frac{L_0}{2},\ y^{-2}N_0\leq N\leq N_0
	\end{align*}
	are seen to satisfy all the above conditions.
	We complete the proof of Theorem 1.1.
	
	\begin{remark}
		For a more rigorous proof, we first fix $Q\leq x^{1/10}$ and $R\leq x^{5/9-O(\varepsilon)}Q^{-14/9}$, and then choose suitable $M,N,L$ subject to all conditions. We may assume $QR>x^{1/2-\varepsilon}$ by Lemma 4.2 below.
	\end{remark}
	
	\section{Proof of Theorem 1.4}
	In our proof, we need the following lemmas.
	\begin{lemma}
		Uniformly for $x\geq y\geq 2$ and $1\leq c<y$, we have
		\begin{align*}
			\Psi(cx,y)=\Psi(x,y)c^{\alpha(x,y)}\Bigg(1+O\left(\frac{1}{u}+\frac{\log y}{y}\right)\Bigg),
		\end{align*}
		where $u=\log x/\log y$,
		\begin{align*}
			\alpha(x,y)=\frac{\log(1+y/\log x)}{\log y}\left(1+O\left(\frac{\log\log(1+y)}{\log y}\right)\right).
		\end{align*}
	\end{lemma}
	\begin{proof}
		This is [11, Corollary 1.7 and (1.26)].
	\end{proof}
	\begin{lemma}
		For any $\varepsilon>0$ and $A>0$, there exists $C>0$ such that for any $y\in[(\log x)^C,x]$, one has
		\begin{align*}
			\sum_{q\leq x^{1/2-\varepsilon}}\max_{(a,q)=1}\Bigg|\Psi(x,y;q,a)-\frac{\Psi_q(x,y)}{\varphi(q)}\Bigg|\ll\frac{\Psi(x,y)}{(\log x)^{A}}.
		\end{align*}
	\end{lemma}
	\begin{proof}
		This is easy from [10, Theorem 1].
	\end{proof}
	
	Consider $(\log x)^C<y<x^{1/C}$ with the $C$ in Theorem 1.1 (or in Lemma 4.2).
	For $n<x/\log x$, we have
	\begin{align*}
		\Psi(x/\log x,y)\asymp\Psi(x,y)(\log x)^{-\alpha(x/\log x,y)}=o(\Psi(x,y)),
	\end{align*}
	since
	\begin{align*}
		\alpha(x/\log x,y)=(1+o(1))\left(1-\frac{\log\log x}{\log y}\right)\geq 1-\frac{1}{C}+o(1).
	\end{align*}
	Considering $x/\log x<n<x$, we thus have
	\begin{align*}
		&\sum_{n\leq x,P^+(n)\leq y}1\\&=\sum_{n\leq x,P^+(n)\leq y}\frac{\log(n+1)}{\log x}+o(\Psi(x,y))\\&=\frac{1}{\log x}\sum_{p'<x}\log p'\sum_{\substack{n\leq x,P^+(n)\leq y\\n\equiv -1\mkern-15mu\pmod{p'}}}1+O\Bigg(\frac{1}{\log x}\sum_{\substack{p',k\\p'^k<x,k\geq 2}}\log p'\sum_{\substack{n\leq x,P^+(n)\leq y\\n\equiv -1\mkern-15mu\pmod{p'^k}}}1\Bigg)+o(\Psi(x,y)).
	\end{align*}
	
	Now we need to estimate the error term. By Lemma 4.2, we have
	\begin{align*}
		&\frac{1}{\log x}\sum_{\substack{p',k\\p'^k<x,k\geq 2}}\log p'\sum_{\substack{n\leq x,P^+(n)\leq y\\n\equiv -1\mkern-15mu\pmod{p'^k}}}1\\&=\frac{1}{\log x}\sum_{\substack{p',k\\p'^k\leq x^{1/3},k\geq 2}}\log p'\sum_{\substack{n\leq x,P^+(n)\leq y\\n\equiv -1\mkern-15mu\pmod{p'^k}}}1+\frac{1}{\log x}\sum_{\substack{p',k\\p'^k>x^{1/3},k\geq 2}}\log p'\sum_{\substack{n\leq x,P^+(n)\leq y\\n\equiv -1\mkern-15mu\pmod{p'^k}}}1\\&\ll \frac{\Psi(x,y)}{\log x}\sum_{\substack{p',k\\p'^k\leq x^{1/3},k\geq 2}}\frac{\log p'}{p'^k}+\frac{x}{\log x}\sum_{2\leq k\leq \log x/\log 2}\sum_{\substack{x^{1/(3k)}<p'<x^{1/2}}}\frac{\log p'}{p'^k}.
	\end{align*}
	For the first summation, we have
	\begin{align*}
		\ll \frac{\Psi(x,y)}{\log x}\sum_{\substack{p'\leq x^{1/3}}}\frac{\log p'}{p'^2}\ll \frac{\Psi(x,y)}{\log x}.
	\end{align*}
	For the second summation, we have
	\begin{align*}
		\ll \frac{x}{\log x}\sum_{2\leq k\leq \log x/\log2}\frac{\log x}{x^{(k-1)/(3k)}}\ll x^{5/6}\log x.
	\end{align*}
	
	Then we have
	\begin{align*}
		\sum_{n\leq x,P^+(n)\leq y}1=\frac{1}{\log x}\sum_{p'<x}\log p'\sum_{\substack{n\leq x,P^+(n)\leq y\\n\equiv -1\mkern-15mu\pmod{p'}}}1+o(\Psi(x,y)).\tag{4.1}
	\end{align*}
	\begin{proposition}
		For $0<\beta<3/5$, we have
		\begin{align*}
			\frac{1}{\log x}\sum_{p'<x^{\beta}}\log p'\sum_{\substack{n\leq x,P^+(n)\leq y\\n\equiv -1\mkern-15mu\pmod{p'}}}1=\left(\beta+o(1)\right)\Psi(x,y).
		\end{align*}
	\end{proposition}
	\begin{proof}
		By Drappeau's work [4, Theorem 1], we have
		\begin{align*}
			&\frac{1}{\log x}\sum_{p'<x^{\beta}}\log p'\sum_{\substack{n\leq x,P^+(n)\leq y\\n\equiv -1\mkern-15mu\pmod{p'}}}1\\&=\frac{1}{\log x}\sum_{p'<x^{\beta}}\frac{\log p'}{p'-1}\sum_{\substack{n\leq x,P^+(n)\leq y\\(n,p')=1}}1+o(\Psi(x,y))\\&=\frac{1}{\log x}\sum_{p'<x^{\beta}}\frac{\log p'}{p'-1}\sum_{\substack{n\leq x,P^+(n)\leq y}}1+O\Bigg(\frac{1}{\log x}\sum_{p'<x^{\beta}}\frac{\log p'}{p'-1}\sum_{\substack{n\leq x,P^+(n)\leq y\\p'|n}}1\Bigg)+o(\Psi(x,y))\\&=\beta\Psi(x,y)+O\Bigg(\frac{1}{\log x}\sum_{p'<y}\frac{\log p'}{p'-1}\Psi(x/p',y)\Bigg)+o(\Psi(x,y)).
		\end{align*}
		For the $O(\cdot)$ term, we have
		\begin{align*}
			&\ll \frac{1}{\log x}\sum_{p'<y}\frac{\log p'}{p'^{1+\alpha(x/p',y)}}\Psi(x,y)\ll\frac{\Psi(x,y)}{\log x},
		\end{align*}
		since
		\begin{align*}
			\alpha(x/p',y)=(1+o(1))\frac{\log y-\log\log(x/p')}{\log y}\geq 1-\frac{1}{C}+o(1).
		\end{align*}
	\end{proof}

	Now we have (The proof is similar to that of (4.1).)
	\begin{align*}
		\sum_{\substack{n\leq x,P^+(n)\leq y}}1&=\sum_{\substack{2\leq n\leq x,P^+(n)\leq y}}\frac{\log(n-1)}{\log x}\frac{\log(n+1)}{\log x}+o(\Psi(x,y))\\&=\frac{1}{(\log x)^2}\sum_{\substack{2\leq n\leq x,P^+(n)\leq y}}\sum_{p_1|(n-1)}\log p_1\sum_{p_2|(n+1)}\log p_2+o(\Psi(x,y)).
	\end{align*}
	For $1/18<\alpha<1/10$, we have
	\begin{align*}
		&\sum_{\substack{2\leq n\leq x,P^+(n)\leq y\\\min(P^+(n-1),P^+(n+1))>x^{1/2+\alpha}}}1\\&\geq \frac{1}{(\log x)^2}\sum_{\substack{2\leq n\leq x,P^+(n)\leq y}}\Bigg(\sum_{p_1|(n-1)}\log p_1\sum_{p_2|(n+1)}\log p_2-\sum_{\substack{p_1|(n-1)\\p_1\leq x^{1/2+\alpha}}}\log p_1\sum_{p_2|(n+1)}\log p_2\\&\quad-\sum_{p_1|(n-1)}\log p_1\sum_{\substack{p_2|(n+1)\\p_2\leq x^{1/2+\alpha}}}\log p_2+\sum_{\substack{p_1|(n-1)\\p_1\leq x^{1/2+\alpha}}}\log p_1\sum_{\substack{p_2|(n+1)\\p_2\leq x^{1/2+\alpha}}}\log p_2\Bigg).
	\end{align*}
	We already have
	\begin{align*}
		\frac{1}{(\log x)^2}\sum_{\substack{2\leq n\leq x,P^+(n)\leq y}}\sum_{p_1|(n-1)}\log p_1\sum_{p_2|(n+1)}\log p_2=(1+o(1))\Psi(x,y),
	\end{align*}
	\begin{align*}
		\frac{1}{(\log x)^2}&\sum_{\substack{2\leq n\leq x,P^+(n)\leq y}}\sum_{\substack{p_1|(n-1)\\p_1\leq x^{1/2+\alpha}}}\log p_1\sum_{p_2|(n+1)}\log p_2\\&=(1+o(1))\frac{1}{\log x}\sum_{p_1\leq x^{1/2+\alpha}}\log p_1\sum_{\substack{n\leq x,P^+(n)\leq y\\n\equiv 1\mkern-15mu\pmod{p_1}}}1=\left(\frac{1}{2}+\alpha+o(1)\right)\Psi(x,y),
	\end{align*}
	\begin{align*}
		\frac{1}{(\log x)^2}\sum_{\substack{2\leq n\leq x,P^+(n)\leq y}}\sum_{\substack{p_1|(n-1)}}\log p_1\sum_{\substack{p_2|(n+1)\\p_2\leq x^{1/2+\alpha}}}\log p_2=\left(\frac{1}{2}+\alpha+o(1)\right)\Psi(x,y).
	\end{align*}
	Then
	\begin{align*}
		&\sum_{\substack{2\leq n\leq x,P^+(n)\leq y\\\min(P^+(n-1),P^+(n+1))>x^{1/2+\alpha}}}1\\&\geq \frac{1}{(\log x)^2}\sum_{\substack{2\leq n\leq x,P^+(n)\leq y}}\sum_{\substack{p_1|(n-1)\\p_1\leq x^{1/2+\alpha}}}\log p_1\sum_{\substack{p_2|(n+1)\\p_2\leq x^{1/2+\alpha}}}\log p_2-(2\alpha+o(1))\Psi(x,y)\\&:=\mathscr{S}_\alpha-(2\alpha+o(1))\Psi(x,y).\tag{4.2}
	\end{align*}
	Changing the order of summation, by Lemma 4.2 and Theorem 1.1, we have
	\begin{align*}
		\mathscr{S}_\alpha\geq &\frac{1}{(\log x)^2}\sum_{\substack{p_1,p_2\\\max(p_1,p_2)\leq x^{1/2+\alpha}\\p_1p_2<x^{1/2-\varepsilon}}}\log p_1\log p_2\sum_{\substack{2\leq n\leq x,P^+(n)\leq y\\n\equiv 1\mkern-15mu\pmod{p_1}\\n\equiv -1\mkern-15mu\pmod{p_2}}}1\\&+\frac{1}{(\log x)^2}\sum_{\substack{p_1,p_2\\\max(p_1,p_2)\leq x^{1/2+\alpha}\\p_1p_2\geq x^{1/2}\\p_1<x^{1/10},p_2<x^{5/9-\varepsilon}p_1^{-14/9}\ \text{or}\ p_2<x^{1/10},p_1<x^{5/9-\varepsilon}p_2^{-14/9}}}\log p_1\log p_2\sum_{\substack{2\leq n\leq x,P^+(n)\leq y\\n\equiv 1\mkern-15mu\pmod{p_1}\\n\equiv -1\mkern-15mu\pmod{p_2}}}1\\=&(1+o(1))\frac{\Psi(x,y)}{(\log x)^2}\Bigg(\sum_{\substack{p_1,p_2\\p_1p_2<x^{1/2-\varepsilon}}}+2\sum_{\substack{p_1,p_2\\\max(p_1,p_2)\leq x^{1/2+\alpha}\\p_1p_2\geq x^{1/2}\\p_1<x^{1/10},p_2<x^{5/9-\varepsilon}p_1^{-14/9}}}\Bigg)\frac{\log p_1\log p_2}{(p_1-1)(p_2-1)}.
	\end{align*}
	We have
	\begin{align*}
		\frac{1}{(\log x)^2}\sum_{\substack{p_1,p_2\\p_1p_2<x^{1/2-\varepsilon}}}\frac{\log p_1\log p_2}{(p_1-1)(p_2-1)}=\frac{1}{8}+O(\varepsilon),
	\end{align*}
	and
	\begin{align*}
		\frac{1}{(\log x)^2}\sum_{\substack{p_1,p_2\\\max(p_1,p_2)\leq x^{1/2+\alpha}\\p_1p_2\geq x^{1/2}\\p_1<x^{1/10},p_2<x^{5/9-\varepsilon}p_1^{-14/9}}}\frac{\log p_1\log p_2}{(p_1-1)(p_2-1)}&=\int_0^{1/10}\left(\frac{1}{18}-\frac{5t}{9}\right)\mathrm{d}t+O(\varepsilon)\\&=\frac{1}{360}+O(\varepsilon).
	\end{align*}
	Finally, we have
	\begin{align*}
		\sum_{\substack{2\leq n\leq x,P^+(n)\leq y\\\min(P^+(n-1),P^+(n+1))>x^{1/2+\alpha}}}1\geq (g(\alpha)+o(1))\Psi(x,y),\tag{4.3}
	\end{align*}
	where $g(\alpha)=1/8+1/180-2\alpha$. It is easy to see that $g(\alpha)>0$ for $\alpha<1/16+1/360$. We complete the proof.
	
	\begin{remark}
		Using our Remark 1.2, a similar argument improves the constant $1/16+1/360\approx 0.0653$ to $0.0661$.
	\end{remark}
	
	\section{Proof of Theorem 1.3}
	
	In this section, our main improvement comes from the following result due to Lichtman [15, Theorem 1.8].
	\begin{lemma}
		Let $C,D,N,R,S\geq 1$. Let $b_{n,r,s}$ be a sequence supported inside $(0,N]\times(R,2R]\times(S,2S]\cap\mathbb{N}^3$. Let $g:\mathbb{R}^5_{>0}\rightarrow \mathbb{C}$ be a smooth function with compactly supported in $(C,2C]\times(D,2D]\times\mathbb{R}^3_{>0}$ such that
		\begin{align*}
			\frac{\partial^{\nu_1+\cdots+\nu_5}}{\partial c^{\nu_1}\partial d^{\nu_2}\partial n^{\nu_3}\partial r^{\nu_4}\partial s^{\nu_5}}g(c,d,n,r,s)\ll_\nu (c^{-\nu_1}d^{-\nu_2}n^{-\nu_3}r^{-\nu_4}s^{-\nu_5})^{1-\varepsilon_0}
		\end{align*}
	   for all $\nu\in (\mathbb{Z}_{\geq 0})^5$ and some small fixed $\varepsilon_0$. Then uniformly for any $a\in \mathbb{Z}_{\neq 0}$, we have
	   \begin{align*}
	   	\sum_{\substack{r\sim R\\s\sim S\\(r,s)=1}}\sum_{0<n\leq N}b_{n,r,s}\sum_{c,d}g(c,d,n,r,s)e\left(an\frac{\overline{rd}}{sc}\right)\ll_{\varepsilon,\varepsilon_0}(|a|CDNRS)^{\varepsilon+O(\varepsilon_0)}\mathscr{J},
	   \end{align*} 
	   where $\|b\|_2=\|\tilde{b}(1)\|_2$, $\|\tilde{b}(n'')\|_2^2=\sum_{n,r,s}|b_{nn'',r,s}|^2$, and $\mathscr{J}=\mathscr{J}(a,C,D,N,R,S)$ is given by
	   \begin{align*}
	   	\mathscr{J}^2=&\sum_{\substack{n''|a^{\infty}\\n''\leq 2N}}(|a|n'')^{\theta}\left(CS\left(\frac{N}{n''}+RS\right)(C+DR)+aNRS\right)\|\tilde{b}(n'')\|_2^2\\&+\left(CS(CD\sqrt{R})^{2\theta}(N+RS)^{1-\theta}(C+DR)^{1-2\theta}+D^2NR\right)\|b\|_2^2.
	   \end{align*}
	\end{lemma}
	\begin{remark}
		The exponent $\theta \leq 7/32$ corresponds to the smallest exceptional eigenvalues $\lambda_1(q) \in (0, 1/4)$ on the Hecke congruence subgroups $\Gamma_0(q)$, $q \geq 1$, by $\theta_q =\sqrt{1/4 - \lambda_1(q)}$. Under Selberg's eigenvalue conjecture there would be no exceptional eigenvalues. In 2003, Kim and Sarnak [14, Appendix 2] showed that $\theta_q \leq 7/32$. (Sometimes, by adjusting the coefficient, this also becomes $7/64$.)
	\end{remark}
	
	Before starting the proof, we collect all the assumptions needed for this section as follows:
	\begin{align*}
		0<a\pp x,\ |\nu d_1d_2|\pp1,\ x^{1-\varepsilon}<MNL\leq x,\ NL\pp K\pp NL,
	\end{align*}
	\begin{align*}
		x^{-20\varepsilon}\ll R\ll x^{-\varepsilon}NL, \ H\ll x^{\varepsilon}RM^{-1}\pp N,\ t\pp 1,\ L^2\pp N^2x,\ N^{3/2}L^{3/2}<x^{1-7\varepsilon}.
	\end{align*}
	(The same result holds when $a$ is negative.) In addition, we have $n,n'\asymp N$, $\ell\asymp L$, $k\asymp K$, $e\asymp E$, $q\asymp Q$.

	We now return to the proof of the theorem.
	Following [4], we only need to estimate (These are the first two equations in [4, p. 847].)
	\begin{align*}
		\mathcal{B}=\sum_e\sum_qB_{e,q}\mathop{\sum_{k}\sum_{\ell}}_{(k,\nu q\ell)=1}\Phi_0\left(\frac{k}{K}\right)\Phi_0\left(\frac{\ell}{L}\right)e\left(ae\frac{\overline{\nu d_2q\ell}}{k}\right),
	\end{align*} 
	where $\Phi_0:\mathbb R\to\mathbb R$ denotes a smooth function bounded above by the indicator function of the interval $[1/2,3]$ and dominating the indicator function of the interval $[1,2]$, the sequence $\beta(n,h)$ satisfies $|\beta(n,h)|\pp 1$,
	\begin{align*}
		B_{e,q}:=\sum_{\substack{n,n'\\nn'=q}}\mathop{\sum\sum}_{\substack{1\leq |h|,|h'|\leq H\\n'h-nh'=e/t}}\beta(n,h)\overline{\beta(n',h')}.
	\end{align*}
	(One small point requires attention: In [4, (3.20) and (3.21)], we may use summation by parts to remove the exponential $e(ah/(\nu k_1r))$, since it is ``flat". We also fix $t$ beforehand.)
	In fact, similar to our discussion in Section 3.3, we only need to prove that (see [4, (3.22)])
	\begin{align*}
		MR^{-1}(KL)^{1/2}\mathcal{B}^{1/2}\ll x^{2-100\varepsilon}(MR)^{-1},
	\end{align*}
	i.e.,
	\begin{align*}
		M^2N^{1/2}L\mathcal{B}^{1/2}\ll x^{2-100\varepsilon}.\tag{5.1}
	\end{align*}
	To simplify the computation, we assume $t \pp 1$, $H \pp NLM^{-1}$.
	The contribution of $e=0$ to $\mathcal{B}$ is
	\begin{align*}
		\mathcal{B}(e=0)\pp KLNH\pp N^3L^3M^{-1}.\tag{5.2}
	\end{align*}
	For $e\neq 0$, we let $n_0 := (n, n')$, substitute $(n_0n, n_0n',n_0q, n_0e) \leftarrow (n, n', \nu d_2q, e)$, and use the triangle inequality to obtain
	\begin{align*}
		\mathcal{B}(e\neq 0)\pp \max_{\substack{N_0\ll N\\E\ll tHNN_0^{-1}\\Q\asymp N^2N_0^{-1}}}\sum_{\substack{q\sim Q}}\max_{\substack{\nu d_2n_0nn'=q\\n,n'\sim N/n_0}}\Bigg|\sum_{e\sim E}B_{e,q}'\mathop{\sum_{k}\sum_{\ell}}_{(k,\nu q\ell)=1}\Phi_0\left(\frac{k}{K}\right)\Phi_0\left(\frac{\ell}{L}\right)e\left(ae\frac{\overline{q\ell}}{k}\right)\Bigg|,
	\end{align*}
	where
	\begin{align*}
		B_{e,q}':=\mathop{\sum\sum}_{\substack{1\leq |h|,|h'|\leq H\\n'(q)h-n(q)h'=e/(tn_0)}}\beta(n(q),h)\overline{\beta(n'(q),h')},
	\end{align*}
	if the maximum is attained at some $n_0(q)$, $n(q)$, $n'(q)$; if the maximum is empty, we let $B_{e,q}':= 0$.
	Applying Lemma 5.1, we have
	\begin{align*}
		\mathcal{B}(e\neq 0)\pp \max_{\substack{N_0\ll N\\E\ll tHNN_0^{-1}\\Q\asymp N^2N_0^{-1}}}\mathscr{J}^2,\tag{5.3}
	\end{align*}
	where $\|B'\|_2=\|\tilde{B}'(1)\|_2$, $\|\tilde{B}'(e'')\|_2^2=\sum_{e,q}|B'_{ee'',q}|^2$, and
	\begin{align*}
		\mathscr{J}^2=&\sum_{\substack{e''|a^{\infty}\\e''\leq 2E}}(ae'')^{\theta}\left(K\left(\frac{E}{e''}+Q\right)(K+LQ)+aEQ\right)\|\tilde{B}'(n'')\|_2^2\\&+\left(K(KL\sqrt{Q})^{2\theta}(E+Q)^{1-\theta}(K+LQ)^{1-2\theta}+L^2EQ\right)\|B'\|_2^2.
	\end{align*}
	For $\|\tilde{B}'(n'')\|_2^2$, we have
	\begin{align*}
		\|\tilde{B}'(n'')\|_2^2\pp& \sum_{n_0\sim N_0}\sum_{n,n'\sim N/n_0}\sum_{\substack{e\sim E\\e''|e}}\Bigg|\sum_{\substack{h,h'\ll H\\n'h-nh'=e/t}}1\Bigg|\\\pp&\sum_{n_0\sim N_0}\sum_{h,h'\ll H}\sum_{\substack{e\sim E\\e''|e}}\sum_{\substack{n\sim N/n_0\\nh'\equiv e/t\mkern-15mu\pmod{h}}}\sum_{n'=(nh'-e/t)/h}\sum_{\substack{h''\ll H\\n'h''\equiv e/t\mkern-15mu\pmod{n}}}1\\\pp&\sum_{n_0\sim N_0}\sum_{h_0\ll H}\sum_{h_1,h_1'\ll H/h_0}\sum_{\substack{e\sim E\\e''|e,h_0|e}}\sum_{\substack{n\sim N/n_0\\nh_1'\equiv e/(th_0)\mkern-15mu\pmod{h_1}}}\left(\frac{N_0H}{N}+1\right)\\\pp&\sum_{n_0\sim N_0}\sum_{h_0\ll H}\sum_{h_1,h_1'\ll H/h_0}\frac{E(e'',h_0)}{e''h_0}\left(\frac{N}{N_0h_1}+1\right)\left(\frac{N_0H}{N}+1\right)\\\pp&e''^{-1}E(N_0H^2+NH+N_0^2H^3N^{-1}).
	\end{align*}
	Since 
	\begin{align*}
		\sum_{e''|a^{\infty}}e''^{\theta-1}\pp1,
	\end{align*}
	we have
	\begin{align*}
		\mathscr{J}^2\pp&\Big\{a^{\theta}\left(K\left(E+Q\right)(K+LQ)+aEQ\right)+\Big(K(KL\sqrt{Q})^{2\theta}(E+Q)^{1-\theta}(K+LQ)^{1-2\theta}\\&+L^2EQ\Big)\Big\}HNN_0^{-1}(N_0H^2+NH+N_0^2H^3N^{-1})\\\pp&\Big\{a^{\theta}\left(K\left(HNN_0^{-1}+N^2N_0^{-1}\right)(K+LN^2N_0^{-1})+aHN^3N_0^{-2}\right)\\&+\Big(K(KLN)^{2\theta}N_0^{-\theta}(HN+N^2)^{1-\theta}N_0^{\theta-1}(K+LN^2/N_0)^{1-2\theta}\\&+L^2HN^3N_0^{-2}\Big)\Big\}HNN_0^{-1}(N_0H^2+NH+N_0^2H^3N^{-1}).
	\end{align*}
	From the monotonicity of $N_0$, by $H\pp \min(R/M,N)$ and $R\pp K\pp NL$, the above becomes (Let $N_0=1$.)
	\begin{align*}
		\mathscr{J}^2\pp&\Big\{a^{\theta}\left(KLN^4+aHN^3\right)\\&+\Big(K(KLN)^{2\theta}N^{2-2\theta}(LN^2)^{1-2\theta}+L^2HN^3\Big)\Big\}HN(H^2+NH+H^3N^{-1})\\\pp&\Big\{a^{\theta}\left(KLN^4+aHN^3\right)+\Big(K(KLN)^{2\theta}N^{2-2\theta}(LN^2)^{1-2\theta}+L^2HN^3\Big)\Big\}H^2N^2\\\pp&\Big\{L^2N^5(x^{\theta}+L^{2\theta}N^{-2\theta})\Big\}N^4L^2M^{-2}\\\pp&M^{-2}N^9L^4(x^{\theta}+L^{2\theta}N^{-2\theta}).\tag{5.4}
	\end{align*}
	Suppose $L^2\pp N^2x$.
	Then by (5.2), (5.3) and (5.4), we have
	\begin{align*}
		M^2N^{1/2}L\mathcal{B}^{1/2}&\pp M^2N^{1/2}L(\mathcal{B}(e=0)+\mathcal{B}(e\neq 0))^{1/2}\\&\pp M^2N^{1/2}L\left((N^3L^3M^{-1})^{1/2}+\left(M^{-2}N^9L^4x^\theta\right)^{1/4}\right)\\&\pp M^{3/2}N^2L^{5/2}+M^{3/2}N^{11/4}L^{2}x^{\theta/4}.
	\end{align*}
	Suppose
	\begin{align*}
		N^{1/2}L\ll x^{1/2-300\varepsilon},\ N^{5/4}L^{1/2}\ll x^{1/2-\theta/4-300\varepsilon},
	\end{align*}
	we have $M^2N^{1/2}L\mathcal{B}^{1/2}\ll x^{2-100\varepsilon}$, which completes our goal in (5.1).
	
	Similar to our discussion in Section 3.3, the best result for $R$ is 
	\begin{align*}
		L\approx x^{(3+\theta)/8},\  R=x^{(5-\theta)/8-310\varepsilon}.
	\end{align*}
	by letting
	\begin{align*}
		M_0:=x^{(3+\theta)/8+300\varepsilon},\ L_0:=x^{(3+\theta)/8-300\varepsilon},\ N_0:=x^{(2-2\theta)/8}.
	\end{align*}
	After verifying all parameters, we complete the proof by letting $\theta=7/32$.

	\section*{Acknowledgements}
	The author thanks Professor Yongxiao Lin for some careful comments. The author also thanks Professor Yongxiao Lin, Jie Wu and Zhiwei Wang for their suggestions on the application of the original version.

\end{document}